\documentclass[reqno,12pt]{amsart} 
\usepackage{vmargin}
\setpapersize{A4}
\usepackage{amsmath} 
\usepackage{amsfonts}   
\usepackage{amssymb}      
\usepackage{eufrak}      
\usepackage{amscd}      
\usepackage{amsthm}      
\usepackage{epsfig}      
\usepackage{amstext}      
\usepackage[all,line,dvips]{xy} 
\CompileMatrices %goer at xy-matricerne koerer hurtigere.

\newcommand{\diag}{\operatorname{diag}}

 \newcommand{\supp}{\operatorname{supp}}

\newcommand{\Dim}{\operatorname{Dim}}
\newcommand{\Rank}{\operatorname{Rank}}

 \newcommand{\tr}{\operatorname{tr}}

   \theoremstyle{plain}%default
   \newtheorem{thm}{Theorem}[section]
   \newtheorem{prop}[thm]{Proposition}
   \newtheorem{lemma}[thm]{Lemma}  
   
   \theoremstyle{definition}
   
   \newtheorem{defn}[thm]{Definition}
   
   \theoremstyle{remark}

\newtheorem{assert}[thm]{Assertion}

\usepackage{graphicx}

   \numberwithin{equation}{section}

        \date{\today}
\title[The crossed product of an AH-algebra  by
  an endomorphism]{Pure infiniteness of the crossed product of an AH-algebra by
  an endomorphism}
  
\author{Klaus Thomsen}

\date{\today}

\email{matkt@imf.au.dk}
\address{Institut for matematiske fag, Ny Munkegade, 8000 Aarhus C, Denmark}

\begin{document}

\maketitle

\section{Introduction}

It has been shown by Deaconu, \cite{De}, and Anantharaman-Delaroche,
\cite{An}, that the
$C^*$-algebra of a local homeomorphism is the crossed product by an
endomorphism of another
$C^*$-algebra. As observed in \cite{De} this implies that such an
algebra is often infinite, and Anantharaman-Delaroche
described in \cite{An} a sufficient condition for the algebra to be purely
infinite. Recall that a simple $C^*$-algebra is said to be purely infinite when all its
non-zero hereditary $C^*$-subalgebras contain an infinite
projection. Thanks to the classification result of Kirchberg and
Phillips this means that the simple and purely infinite $C^*$-algebras which
arise from local homeomorphisms are classified by their K-theory
groups, and it becomes therefore an important question to decide when the
algebra of a local homeomorphism is simple and purely
infinite. 

In \cite{R1} R\o rdam proved that the crossed product by a full corner
endomorphism of a simple unital
$C^*$-algebra of real rank zero with comparability of projections is simple and purely infinite. In particular,
the crossed product of a simple unital AF-algebra by such an endomorphism is
simple and purely infinite. In the same paper R\o rdam initiated also the
classification of purely infinite simple $C^*$-algebras which was
subsequently completed, mutatis mutandis, by the classification results of Kirchberg and Phillips mentioned above. R\o rdams result on the
crossed product by an endomorphism has been extended and used by several
other mathematicians, but in most of these results the initial
algebra, the one with the endomorphism, hereafter called the
\emph{core}, is assumed to be simple and to have various other properties. The simplicity of the crossed product, as well as its pure
infiniteness, is then a consequence. The work of Dykema and R\o rdam
in \cite{KR} is an exception, but they assume some rather special
properties of the endomorphism which are not easy to establish.

For the application to the $C^*$-algebras of a local homeomorphism it
is a nuisance to have to assume simplicity of the core. When the algebra of a local homeomorphism is simple, the core may or may not be simple
and hence the existing results, as the one of R\o rdam, on crossed
products by endomorphisms can generally only be used by imposing
additional assumptions. It is the
purpose of the present paper to obtain a result about the pure
infiniteness of a crossed product by an endomorphism in which
simplicity is assumed of the crossed product rather than of the core,
and which is general enough to cover the $C^*$-algebra of a local
homeomorphism; assuming only that it is simple. The following, which
is the main result of the paper, is such a theorem. The definition of
a 'unital AH-algebra with slow
  dimension growth' will be given in the next section.

\begin{thm}\label{semicrossed} Let $A$ be a unital AH-algebra with slow
  dimension growth. Let $\beta : A \to A$ be an injective endomorphism such that
\begin{enumerate}
\item[i)] $\beta(1)$ is a full projection in $A$
  (i.e. $\overline{A\beta(1)A} = A$), and
\item[ii)] there is no trace state $\omega$ of $A$ such that $\omega
  \circ \beta = \omega$.%, and 
%\item[iii)] $A$ is tracially almost divisible.
\end{enumerate}
If the crossed product $A \times_{\beta} \mathbb N$ is simple, it is also purely infinite.
\end{thm}

Applications of this result to the $C^*$-algebras of local homeomorphisms and locally
injective surjections will be given in \cite{CT}.

It must be observed that the crossed product $A \times_{\beta} \mathbb N$
in the theorem is not the same as the one which was introduced by
W. Paschke and used by R\o rdam in
\cite{R1} where it is assumed that
$\beta$ maps \emph{onto} the corner $\beta(1)A\beta(1)$. In order to cover
also the crossed products by endomorphisms arising from a locally
injective surjection, which may not be open and hence not a local
homeomorphism, cf. \cite{Th1}, we use instead the crossed product
introduced by Stacey in \cite{St}. It can be defined as the universal $C^*$-algebra
generated by a copy of $A$ and an isometry $v$ with the property that
$vav^* = \beta(a)$, \cite{BKR}, and hence it agrees with the one used
by Dykema and R\o rdam in \cite{DR}. Compared to the crossed product of
Paschke, it is not required that $v^*Av \subseteq
A$. When $\beta$ maps onto $\beta(1)A\beta(1)$, as is for example the
case when the situation arises from a local homeomorphism as in
\cite{An}, the two crossed products coincide.

The main strategy of the proof is due to R\o rdam. In \cite{R2} he
proved that the crossed product of a $C^*$-algebra $A$ by an
automorphism is (simple and) purely infinite when $A$ is
\begin{enumerate}
\item{-} exact, finite and separable,
\item{-} simple, 
\item{-} approximately divisible,
\end{enumerate}   
and has no densely defined non-zero trace which is invariant under
the given automorphism. Although this is a result about an
automorphism it has bearing on crossed products by endomorphisms since
they can be realised as a corner in a crossed product by an
automorphism. 

The last condition in the above statement, about the absence of invariant
traces, is of course necessary. The first conditions (1) are harmless and satisfied when the crossed
product arises from one of the locally injective surjections we have in
mind. As we explained above the assumed simplicity is an assumption we aim to
move from the core to the crossed product, while approximate divisibility
is a property which is hard to establish and about which we know next
to nothing when the algebra comes from a local homeomorphism and the
core is not simple. It is
therefore interesting to observe that an important step in the following
proof of Theorem \ref{semicrossed} will be to show that a much weaker
version of divisibility is automatic for unital AH-algebras with slow
dimension growth.

\section{Tracial almost divisibility for AH-algebras with slow
  dimension growth}

Let $M_l$ denote the $C^*$-algebra of complex $l \times l$-matrices.
In the following a \emph{homogeneous $C^*$-algebra} will be a
$C^*$-algebra $A$ isomorphic to a $C^*$-algebra of the form
$$
eC(X,M_l)e
$$
where $X$ is a compact metric space and $e$ is a projection in
$C(X,M_l)$ such that $e(x) \neq 0$ for all $x \in X$. The \emph{dimension ratio} $r(A)$ of $A$ is then defined
to be the number
$$
r(A) = \max_{x \in X} \frac{\Dim X +1}{\Rank e(x)}.
$$

\begin{defn}\label{slowdim} A unital $C^*$-algebra $A$ is an \emph{AH-algebra} when
  there is an increasing sequence $A_1 \subseteq A_2 \subseteq A_3
  \subseteq \dots$ of unital $C^*$-subalgebras of $A$ such that $A =
  \overline{\bigcup_n A_n}$ and each $A_n$ is a homogeneous $C^*$-algebra. We say that $A$ has \emph{slow dimension
    growth} when there is such a sequence with the additional property
  that $\lim_{n \to \infty} r(A_n) = 0$. 
\end{defn}
There seems to be slightly varying definitions of slow dimension growth
for AH-algebras and it should therefore be observed that with the above
definition we insist that the rank of the projections increase without
bounds even when all the involved topological spaces are
zero-dimensional.

%Even though we only aim to deal with endomorphisms of AH-algebras with
%slow dimension growth, for the proof we shall need the following
%'local' version.

%\begin{defn} A unital $C^*$-algebra $A$ is \emph{almost an AH-algebra with
 %   slow dimension growth} when the following holds: For any $x \in A$
 % and any $\epsilon > 0$ there is a unital homogeneous
 % $C^*$-subalgebra $1 \in B \subseteq A$ and an element $y \in B$ such
 % that $r(B) \leq \epsilon$ and $\left\| x- y\right\| \leq \epsilon$.
%\end{defn}

%Of course, a unital AH algebra with slow dimension growth is also
%almost so, but the converse is probably false.

Let $A$ be a $C^*$-algebra and $a,b$ two positive elements of
$A$. Recall, cf. e.g. \cite{T}, that $a$ is \emph{Cuntz subequivalent}
to $b$ when there is a sequence $\{z_n\}$ in $A$ such that $a =
\lim_{n \to \infty} z_nbz_n^*$. We write $a \preceq b$ when this
holds. This notion extends the well known subequivalence
in the sense of Murray-von Neumann used for projections. 

In the following we denote by $T(A)$ the convex set of trace states of
a unital $C^*$-algebra $A$. The next definition is inspired by
Definition 2.5 (ii) of \cite{W}.

\begin{defn}\label{tracedivi} A unital $C^*$-algebra $A$ is
  \emph{tracially almost divisible} when the following holds: For any
  positive contraction $h$ in $A$ and any given $m \in \mathbb N$
  there is a $\delta > 0$ with the property that for all $\epsilon > 0$
  there are mutually orthogonal positive contractions
  $h_1,h_2, \dots, h_m$ in $A$ such that 
$$
h_1+h_2 + \dots + h_m \preceq
  h
$$ and 
$$
\tau(h_i) \geq \delta\tau(h) - \epsilon
$$ 
for all $i$
  and all $\tau \in T(A)$.
\end{defn}

As an important step towards the main result of the paper we prove
first the following.

%begin{prop}\label{tracAH} Let $A$ be a unital $C^*$-algebra. Assume $A$ is
%  almost an AH-algebra with slow dimension growth. Then $A$ is tracially almost divisible.
%\end{prop}

\begin{prop}\label{tracAH} Let $A$ be a unital AH-algebra with slow
  dimension growth. Then $A$ is tracially almost divisible.
\end{prop}
We will actually prove a slightly stronger result; namely that the $\delta$ of Definition \ref{tracedivi}
can be chosen to be $\frac{1}{4m}$, independently of $h$. However, the
proof of the main result will not require this strengthening of the
conclusion.

The main tools for the proof of Proposition \ref{tracAH} are methods
and results of
A. Toms from \cite{T} about Cuntz subequivalence in a homogeneous
$C^*$-algebra. When $A$ is a unital $C^*$-algebra and $\tau \in T(A)$
there is associated to $\tau$ a 'dimension function' $d_{\tau}$ defined on positive
contractions of $A$ as
$$
d_{\tau}(a) = \lim_{n \to \infty} \tau\left(a^{\frac{1}{n}}\right) =
  \sup_{n \in \mathbb N} \tau\left(a^{\frac{1}{n}}\right) .
$$
By Corollary 5.2 of \cite{T} we have the following

\begin{thm}\label{TOMS}(A. Toms) Let $A \simeq eC(X,M_l)e$ be a homogeneous
  $C^*$-algebra. Let $a,b \in A$ be positive contractions such that 
$$
d_{\tau}(a) + \max_{x \in X} \frac{\Dim X}{2\Rank e(x)}  \leq
d_{\tau}(b)
$$
for all $\tau \in T(A)$. It follows that $a \preceq b$.
\end{thm}

Actually, the result in \cite{T} is slightly stronger, but the above
theorem suffices for our purposes.

\begin{lemma}\label{anvtoms}
Let $\epsilon \in ]0,\frac{1}{8}[$ and $m \in \mathbb N$. Let $eC(X,M_l)e$ be a
homogeneous $C^*$-algebra such that $\Rank e(x) = M$ is constant and
$$
\frac{\Dim X +1}{M} < \frac{\epsilon}{8m}.
$$
It follows that for every positive contraction $h \in eC(X,M_l)e$
there are $m$ mutually orthogonal positive contractions $h_1,h_2,
\dots, h_m$ in $eC(X,M_l)e$ such that 
$$
h_1 + h_2 + \dots + h_m \preceq
h
$$ 
and 
\begin{equation}\label{tracelo}
\tau(h_i) \geq \frac{1}{4m}\tau(h) - 2 \epsilon
\end{equation} 
for all $i$
and all trace states $\tau \in T(eC(X,M_l)e)$.   
\end{lemma}
\begin{proof} Let $ j \in \mathbb N, \ j \geq 2$, and set $d = \Dim X +1$. Since $\frac{d}{M} <
  \frac{\epsilon}{4}$ we find that
\begin{equation*}
\begin{split}
& \frac{(j -
  \frac{1}{2})\epsilon}{m} - \frac{d}{2Mm} - \frac{(j +
  \frac{1}{2})\epsilon}{2m} > \frac{(j -
  \frac{1}{2})\epsilon}{m} - \frac{\epsilon}{8m} - \frac{(j +
  \frac{1}{2})\epsilon}{2m} = \frac{(j
  -\frac{7}{4})\epsilon}{2m} \geq \frac{\epsilon}{8m}. 
\end{split}
\end{equation*}  
As $\frac{1}{M} <  \frac{\epsilon}{8m}$ these estimates show that there is a natural number $\alpha_j$
such that
\begin{equation*}
\frac{(j +
  \frac{1}{2})\epsilon}{2} \leq\frac{m\alpha_j}{M}   < (j -
  \frac{1}{2})\epsilon  - \frac{d}{2M}.
\end{equation*}
Note that we can arrange that $\alpha_j \leq \alpha_{j+1}$. Let $J$ be
the least natural
number such that $\left(J+\frac{1}{2}\right)\epsilon \geq
\frac{1}{2}$. (The condition $\epsilon < \frac{1}{8}$ ensures that $J
\geq 4$.) Then 
$$
\frac{m\alpha_j}{M} \leq \frac{m\alpha_J}{M} \leq (J -
  \frac{1}{2})\epsilon  - \frac{d}{2M} < 1
$$
for $2 \leq j \leq J$. We can
therefore choose mutually orthogonal trivial projections $p^j_1,p^j_2,
\dots, p^j_m$ in $C(X,M_l)$ for each $2 \leq j \leq J$ such that $\Rank p^j_i = \alpha_j, i = 1,2,\dots,m$, and such that 
$$
p^j_i \leq p^{j+1}_i, \ i = 1,2, \dots, m, \ 2 \leq j \leq J-1 .
$$ 
Then 
$$
\Rank \left(\sum_{i=1}^m p^J_i\right) + \frac{d}{2}=
{m\alpha_J} + \frac{d}{2} \leq (J -
  \frac{1}{2})\epsilon M \leq \frac{M}{2} \leq M  = \Rank e
$$
so the projection $\sum_{i=1}^m p^J_i$ is Murray-von Neumann
equivalent to a subprojection of $e$; this is a classical fact about vector
bundles, but it follows also from Theorem \ref{TOMS}. Consequently we may assume that $p_i^j
\in eC(X,M_l)e$ for all $i,j$, with the reduction that they may no
longer be trivial projections. Set $p^0_i = p^1_i = 0$ for $i = 1,2,
\dots,m$.

For each $0 \leq j \leq J-1$ we choose a continuous function $g_j :
\left[(j+\frac{1}{2})\epsilon, (j + \frac{3}{2})\epsilon\right] \to [0,1]$
  such that $g_j((j+\frac{1}{2})\epsilon) = 1$ and
  $g_j\left((j+\frac{3}{2})\epsilon\right) = 0$. For each $i =
  1,2,\dots, m$, define a continuous function
$$
H_i : [0,1] \times X \to M_l
$$
such that 
$$
H_i(t,x) = \begin{cases} 0, & \ t \in \left[0,\frac{1}{2}\epsilon\right], \\
  g_j(t)p^j_i(x) + (1-g_j(t)) p^{j+1}_i(x), & \ t \in
  \left[(j+\frac{1}{2})\epsilon, (j + \frac{3}{2})\epsilon\right] , \ 0
  \leq j \leq J-1, \\ p^J_i(x), & \ t \geq (J + \frac{1}{2})\epsilon . 
\end{cases}   
$$
Then $H_i(t,x)$ is a positive contraction and $H_i(t,x)H_{i'}(t,x) =
0, \ i \neq i'$, for all $t,x$.

For each $x \in X$ we consider the extremal trace state $\tau_x$ of
$eC(X,M_l)e$ defined as $\tau_x(f) = \tr (f(x))$ where $\tr$ is the
trace state of $e(x)M_le(x)$. For each $i = 1,2,\dots, m$ we define $h'_i \in C(X,M_l)$ such that
$$
h'_i(x) = H_i\left(\tau_x(h),x \right) . 
$$
Note that the $h'_i$'s are mutually orthogonal positive contractions
and that $h'_i \in eC(X,M_l)e$ since $e(x)H_i(t,x)e(x) = H_i(t,x)$ for all $t,x$.

%We use next the $p^j_i$'s to construct orthogonal contractions $\Phi_i
%\in eC(X,M_l)e$. 
%For each $x \in X$ let $\tau_x \in T(eC(X,M_l)e)$ be the trace defined
%such that $\tau_x(a) = \tr (a(x))$ where $\tr$ is the trace state of
%$e(x)M_le(x)$. Let $j_0 \in \mathbb N$ be the least natural number
%such that $\left\{ x\in X : \tau_x(h) \leq \left(j_0 +
%    \frac{1}{2}\right)\epsilon \right\} \neq \emptyset$. When
%$\tau_x(h) \leq \left(j_0 + \frac{1}{2}\right)$, set $\Phi_i(x) =
%p^{j_0}_i$. For each $j_0 \leq j  \leq J-1$ choose a continuous
%function $g_j : X \to [0,1]$ such that $g_j(x) = 1$ when $\tau_x(h)
%\leq \left(j+ \frac{1}{2}\right)\epsilon$ and $g_j(x) = 0$ when
%$\tau_x(h) > \left(j + \frac{3}{2}\right)\epsilon$. Set
%$$
%\Phi_i(x) = g_j(x)p^j_i + (1-g_j(x))p^{j+1}_i
%$$
%when $\tau_x(h) \in \left[\left(j+\frac{1}{2}\right)\epsilon, \left(j
%    + \frac{3}{2}\right)\epsilon\right]$. When $\tau_x(h) \geq \left(J
%  + \frac{1}{2}\right)\epsilon$ we set $\Phi_i(x) = p^J_i$. Note that
%we assumed that $j_0 \leq J-1$ in this construction. Should this fail
%we set $\Phi_i(x) = p^J_i$ for all $x$. 

Set 
$$
U = \left\{ x \in X : \ \tau_x(h) > (2 + \frac{1}{2})\epsilon \right\}
$$
and consider an $x \in \overline{U}$. Then $\tau_x(h) \in
\left[(j+\frac{1}{2})\epsilon, (j + \frac{3}{2})\epsilon\right]$ for
some $2 \leq j \leq J-1$ or $\tau_x(h) \geq
(J+\frac{1}{2})\epsilon$. In the first case we find that
\begin{equation*}
\begin{split}
&\tau_x(h'_i(x)) = g_j\left(\tau_x(h)\right)\frac{\alpha_j}{M}  +
  \left( 1 - g_j\left(\tau_x(h)\right)\right)  \frac{\alpha_{j+1}}{M}
  \\
&\geq  g_j\left(\tau_x(h)\right)\frac{(j+\frac{1}{2})\epsilon}{2m}  +
  \left( 1 - g_j\left(\tau_x(h)\right)\right)
  \frac{(j+\frac{3}{2})\epsilon}{2m} \geq \frac{\tau_x(h)}{2m} - \frac{\epsilon}{2m},
\end{split}
\end{equation*}
while
$$
d_{\tau_x}\left(\sum_{i=1}^m h'_i\right) \leq \frac{m \alpha_{j+1}}{M} \leq (j
  + \frac{1}{2})\epsilon - \frac{d}{2M} \leq \tau_x(h) - \frac{d}{2M} .
$$
When $\tau_x(h) \geq (J + \frac{1}{2})\epsilon $ we find that
$$
\tau_x(h'_i) = \frac{\alpha_J}{M} \geq \frac{(J +
  \frac{1}{2})\epsilon}{2m} \geq \frac{1}{4m} \geq
\frac{1}{4m}\tau_x(h)
$$
 while 
$$
d_{\tau_x}\left(\sum_{i=1}^m h'_i\right)  
 \leq \frac{m \alpha_J}{M}  \leq
(J-\frac{1}{2}) \epsilon - \frac{d}{2M} \leq \tau_x(h) - \frac{d}{2M}.
$$
All in all we conclude that
$$
\tau_x(h'_i) \geq \frac{1}{4m}\tau_x(h) - \frac{\epsilon}{2}
$$
and
$$
d_{\tau_x}\left(\sum_{i=1}^m h'_i\right) \leq \tau_x(h) -
\frac{d}{2M}
$$
for all $i = 1,2,\dots,m$ and all $x \in \overline{U}$.

If $U= \emptyset$, we
set $h_1= h_2 = \dots = h_m = 0$. Since $0 \preceq h$ and
$\frac{1}{4m}\tau_x(h) \leq 2\epsilon$ for all $x$ this will prove the
lemma in this case. Assume therefore that $U \neq \emptyset$. Recall
that $T\left(eC(\overline{U},M_l)e\right)$ is the closed convex hull of
$\left\{\tau_x : x \in \overline{U} \right\}$. Since $d_{\tau_x} \left(\sum_{i=1}^m h'_i\right) \leq \tau_x(h) - \frac{d}{2M}$
for all $x \in \overline{U}$, and since 
$$
\tau \mapsto
d_{\tau}\left(\sum_{i=1}^m h'_i|_{\overline{U}}\right) - \tau(h|_{\overline{U}})+ \frac{d}{2M}
$$ 
is affine and lower
semi-continuous on $T\left(eC(\overline{U},M_l)e\right)$ we find that 
$$
d_{\tau}\left(\sum_{i=1}^m h'_i|_{\overline{U}}\right) + \frac{d}{2M} \leq \tau(h|_{\overline{U}})  \leq d_{\tau}(h|_{\overline{U}}) 
$$
for all $\tau \in T\left(eC(\overline{U},M_l)e\right)$. Since $\Dim
\overline{U} \leq \Dim X$, it follows from Theorem \ref{TOMS} that there is a
sequence $\{z_n\} \in eC(\overline{U},M_l)e$ such that
$$
\lim_{n \to \infty} z_n(x)h(x)z_n(x)^* = \sum_{i=1}^m h'_i(x)
$$
uniformly in $x \in \overline{U}$. Set
$$
K = \left\{ x \in X : \ \tau_x(h) \geq (3+\frac{1}{2})\epsilon 
\right\}
$$ 
and let $\psi : X
\to [0,1]$ be a continuous function such that $\psi(x) = 1, x \in K$,
and $\supp \psi \subseteq U$. We consider $\psi$ as a central element
of $eC(X,M_l)e$ in the obvious way. Let now
$$
h_i = \psi h'_i, \ i = 1,2, \dots, m,
$$ 
and set $z'_n= z_n\sqrt{\psi} \in eC(X,M_l)e$. Then 
$$
\lim_{n \to \infty} z'_nh{z'_n}^* = \sum_{i=1}^m
h_i
$$ 
and $\tau_x(h_i) \geq \frac{1}{4m}\tau_x(h) - \epsilon$ for
all $x \in K$. Since$\frac{1}{4m}\tau_x(h) - 2\epsilon \leq
\frac{4\epsilon}{4m} - 2\epsilon < 0 \leq \tau_x(h_i)$ when $x \notin K$, we obtain (\ref{tracelo}).

\end{proof}

\emph{Proof of Proposition \ref{tracAH}}: Consider a positive
  contraction $b \in A$. Let $m \in \mathbb N$ and $\epsilon >
  0$ be given. We will complete the proof by showing that there
  are mutually orthogonal positive contractions $b_1,b_2, \dots, b_m$
  in $A$ such that $b_1 + b_2 + \dots + b_m \preceq b$ and $\tau(b_i)
  \geq \frac{1}{4m} \tau(b) - 3 \epsilon$ for all $i$.

Since $A$ is an AH-algebra with slow dimension growth it follows from Lemma 2.5 (ii) in \cite{KR}
  that there is a unital homogeneous $C^*$-sub-algebra $1 \in B
  \subseteq A$ and a positive contraction $a
  \in B$ such that $r(B) < \frac{\epsilon}{8m}$, $a \preceq b$ and $\left\| a-b\right\|\leq
  \epsilon$. Note that $B$ is isomorphic to a direct sum
$$
B \simeq \oplus_{j=1}^N e_jC(X_j,M_l)e_j
$$
of homogeneous $C^*$-algebras such that $\Rank e_j(x) = K_j$ is
constant on $X_j$ and 
$$
\frac{\Dim X_j +1}{K_j} \leq r(B) < \frac{\epsilon}{8m}
$$
for all $j$. We can therefore apply Lemma \ref{anvtoms} to each
summand and in this way obtain mutually orthogonal positive contractions $b_1,b_2, \dots, b_m$
  in $B$ such that $b_1 + b_2 + \dots + b_m \preceq a$ and
  $\tau(b_i) \geq \frac{1}{4m}\tau(a) - 2\epsilon$ for all $\tau \in T(B)$. Since
  $\frac{1}{4m}\tau(a) \geq \frac{1}{4m}\left(\tau(b) -
    \epsilon\right) \geq \frac{1}{4m} \tau(b) - \epsilon$ for all
  $\tau \in T(A)$, we are done.    \qed

\section{Proof of the main result}

In this section we prove Theorem \ref{semicrossed} by an elaboration
of R\o rdams proof of Theorem 2.1 in
\cite{R2}. For this purpose we isolate the following lemmas. In the statement of the
first we use
the (standard) notation $M_{\infty}(B)$ for the union $\bigcup_n
M_n(B)$. Recall that a projection $p \in B$ is
\emph{full} when $\overline{BpB} = B$.

\begin{lemma}\label{tomsvector} Let $Z$ be a compact metric space of
  dimension  $\dim Z \leq d$ and let
  $p \leq e \in C(Z,M_l)$ be projections. Assume that there is a
  natural number $N \in \mathbb N$ such that 
$$
\Rank p(z) > (N+1)(N+2)\left[\frac{d}{2}\right]
$$
for all $z \in Z$, where $[\frac{d}{2}]$ is the least
natural number larger or equal to $\frac{d}{2}$.  It follows that there is a projection
$p' \in M_{\infty}\left(eC(Z,M_l)e\right)$ such that
\begin{equation}\label{tomschok1}
N[p'] \leq [p] \leq (N+3)[p']
\end{equation}
in $K_0(eC(Z,M_l)e)$

 \end{lemma}
\begin{proof}  Let $Z = Z_1 \sqcup Z_2 \sqcup \dots \sqcup Z_k$ be a
  partition of $Z$ by clopen sets such that $\Rank p$ is constant on
  each $Z_j$. Fix $j$ and set $d_j = \Dim Z_j$. Note that $d_j \leq d$. Write $\Rank p = l(N+1)\left[\frac{d_j}{2}\right]
  + r$ where $l,r \in \mathbb N$ and $1 \leq r \leq
  (N+1)\left[\frac{d_j}{2}\right]$. Then $l \geq N+2$ by assumption and hence 
\begin{equation}\label{tomschok}
Nl\left[\frac{d_j}{2}\right] + l\left[\frac{d_j}{2}\right] < \Rank p \leq Nl\left[\frac{d_j}{2}\right] + 2l\left[\frac{d_j}{2}\right]
\end{equation}
on $Z_j$. Let $q_j$ be a trivial projection on $Z_j$ of constant
rank $l\left[\frac{d_j}{2}\right]$. Since $e|_{Z_j}$ is a full projection, $q_j$
is equivalent to a projection $p'_j$ in $M_d\left(eC(Z_j,M_l)e\right)$ for
some $d$. Since $l\left[\frac{d_j}{2}\right]
\geq \frac{d_j}{2}$ it follows from (\ref{tomschok}) and the theory of
vector bundles (or
Theorem \ref{TOMS}), that $N[p'_j] \leq \left[p|_{Z_j}\right] \leq
(N+3)[p'_j]$ in $K_0(eC(Z_j,M_l)e)$. Set $p' =
\sum_j p'_j$. 
\end{proof}

\begin{lemma}\label{typisk} Let $A$ be a unital AH-algebra with slow
  dimension growth and $A_1 \subseteq A_2 \subseteq A_3 \subseteq
  \dots $ a sequence of homogeneous $C^*$-sub-algebras,
$$
A_n \simeq e_nC\left(X_n,M_{m_n}\right)e_n,
$$
such that $1 \in A_1$, $A = \overline{\bigcup_n A_n}$ and 
$
\lim_{n \to
  \infty} r(A_n) = 0$. Let $p$ be a full projection in $A$ and let $K \in \mathbb N$ be given. It
follows that there is an $n \in \mathbb N$ and a projection $q \in
A_n$ such that $q$ is unitarily equivalent to $p$ in $A$ and 
$$
\Rank q(x) \geq K \left(\Dim X_n +1\right) 
$$
for all $x \in X_n$.
\begin{proof} A standard argument shows that $p$ is unitarily
  equivalent to a projection $q$ in $A_m$ for some $m$. Since $p$ is full
  we can assume, by increasing $m$, that $q$ is full in $A_m$. There
  is then a $k \in \mathbb N$ such that $\Rank q(x) \geq
  \frac{1}{k}\Rank e_n(x)$ for all $x \in X_n$, $n \geq m$. Therefore
  the desired inequality will hold for all sufficiently large $n$
  thanks to the slow dimension growth condition.  
\end{proof}
\end{lemma}

Let $a,b \in K_0(A)$. In the following we write $a
\prec b$ when there is a full projection $q$ in $M_{n}(A)$ for some
$n$ such
that $b-a = [q]$.  

\begin{lemma}\label{divisibility} Let $A$ be an AH-algebra with slow
  dimension growth. Let $e,f \in A$ be projections such that $N[e]
  \prec N[f]$ in $K_0(A)$ for some $N \in \mathbb N$. It follows that $[e]  \prec [f]$ in
  $K_0(A)$ and $e \preceq f$ in $A$.
 \end{lemma}
\begin{proof} It follows from Lemma \ref{typisk} that the difference
  $N[f] - N[e]$ is represented by a projection $p$ in a homogeneous
  $C^*$-algebra, containing also projections $e'$ and $f'$, unitarily
  equivalent to $e$ and $f$, respectively, such
  that $\inf_x \Rank p(x)$ is greater than $N +1$ times the dimension 
  of the spectrum. Then well known facts about vector
  bundles, or Theorem \ref{TOMS}, show
  that in this algebra $e'$ is equivalent to a subprojection $e''$ of
  $f'$ such that $f'-e''$ is full. The lemma follows.
\end {proof}

\begin{lemma}\label{fullinf} Let $B$ be a $C^*$-algebra with the
  property that $B = \overline{\bigcup_n B_n}$ where $B_1 \subseteq
  B_2 \subseteq B_3 \subseteq \dots$ are
  $C^*$-subalgebras of $B$ each of which is a unital AH-algebra
  with slow dimension growth. Furthermore, assume that the unit of
  $B_n$ is a full projection in $B_{n+1}$ for each $n$. Let $\alpha$ be an automorphism of $B$
  such that $\omega \circ \alpha \neq \omega$ for all non-zero densely
  defined lower semi-continuous traces $\omega$ on $B$. Assume that $B \times_{\alpha}
  \mathbb Z$ is simple. It follows that every full projection
  in $B$ is infinite in $B \times_{\alpha} \mathbb Z$.
\end{lemma}
\begin{proof} We elaborate on R\o rdams proof of Lemma 2.5 in
  \cite{R2}. Let $e$ be a full projection in $B$. 

\smallskip
a) The first step is to
  show that there is an element $x \in K_0(B)$ such that $x \geq
  \alpha_*(x)$ and $x \neq \alpha_*(x)$. As in \cite{R2} this follows from the absence of
  $\alpha$-invariant traces, by use of results of Blackadar, R\o rdam
  and Goodearl, Handelman. We refer to \cite{R2} for the details of
  the argument.

\smallskip

b) Let $I \subseteq B$ be a non-zero ideal such that $\alpha(I)
\subseteq I$. It follows that $I = B$. Indeed 
$$
J =\overline{\bigcup_{n \geq 0} \alpha^{-n}(I)}
$$
is an ideal in $B$ such that $\alpha(J) = J$. Since $B \times_{\alpha}
\mathbb Z$ is simple it follows that $J = B$. In particular there is an $n \in \mathbb N$ and an element $b
\in I$ such that $\left\|\alpha^n(e) - b \right\| = \left\|e -
  \alpha^{-n}(b)\right\|  < \frac{1}{3}$. As is well-known this
implies that $I$ contains a projection equivalent to
$\alpha^n(e)$. This projection is full in $A$ since $\alpha^n(e)$
is, whence $I = B$.

%(An example: Let $A = C_0(\mathbb Z)$ and $\alpha$ the left-shift,
%$\alpha(f)(z) = f(z+1)$. Then $I = \left\{ f \in A : \ f(z) = 0, \ z
%  \geq 0\right\}$ is a non-trivial ideal such that $\alpha(I)
%\subseteq I$. Note that $C_0(\mathbb Z) \times_{\alpha} \mathbb Z
%\simeq \mathbb K$ is simple and that $C_0(\mathbb Z)$ does not contain
%any full projection!)

\smallskip

In the following we extend $\alpha$ to $M_n(B)$ for all $n$ in the
canonical way. 

c) Let $p \in M_n(B)$ be a projection such that $[p] = x -
\alpha_*(x)$ where $x\in K_0(B)$ is the element from a). Since
$$
\overline{\bigcup_k  \left\{ a_0pb_0 + a_1\alpha(p)b_1 + \dots +
    a_k\alpha^k(p)b_k: \ a_i,b_i \in M_n(B), \ i = 0,1, \dots, k \right\} }
$$
is a non-zero ideal $I$ in $M_n(B)$ such that $\alpha(I) \subseteq I$, it follows
from b) that $I = M_n(B)$ and hence it contains a full projection. By
definition of $I$ this implies that there is a $k$ such that $[p] + \alpha_*[p] +
\alpha^2_* [p] + \dots + \alpha^k_*[p]$ is an order-unit in
$K_0(B)$. Set $y = x + \alpha_*(x) + \alpha_*^2(x) + \dots +
\alpha_*^k(x)$. Then $y -
\alpha_*(y) = [p] + \alpha_*[p] +
\alpha^2_* [p] + \dots + \alpha^k_*[p]$. By exchanging $y$ for $x$ we
may therefore assume that $x - \alpha_*(x) = [p]$ for some full projection $p$
of $M_{\infty}(B)$, i.e. $\alpha_*(x) \prec x$.

\smallskip

d) Write $x = g_1 - g_2$ where $g_1,g_2 \in K_0(B)^+$. We may assume
that $g_i \succ 0, i =1,2$. There is an $N
\in \mathbb N$ such that
$$
3g_1 + 3\alpha_*(g_2) \prec N(x - \alpha_*(x))
$$
and
$$
g_1 + \alpha_*(g_2) \prec N[e] .
$$
Note that since $g_i \geq 0$ there are $L,n \in \mathbb N$ such that
$g_i = [p_i]$ for some projection $p_i$ in $M_L(B_n)$. Since in fact
$g_i \succ 0$ it follows from the assumption about the unit of $B_n$ being
full in $B_{n+1}$ for each $n$, that we can assume that $p_i$ is full
in $M_L(B_n)$. It follows then from Lemma
\ref{typisk} that we can realise the projections $p_i, i = 1,2$, in a homogeneous
$C^*$-subalgebra of $M_L(B_n)$ such that the assumptions of Lemma
\ref{tomsvector} hold for both. In
this way we get elements $f_1,f_2 \in K_0(B)^+$ such that $Nf_j \leq g_j
\leq (N+3)f_j, j = 1,2$, and we set $f= f_1 -f_2$. Then
$$
N[p] \succ  3g_1 + 3\alpha_*(g_2) \geq 3N\left(f_1 +
  \alpha_*(f_2)\right)
$$
which implies that $[p] \succ 3\left(f_1 +
  \alpha_*(f_2)\right)$ by Lemma \ref{divisibility}.
It follows first that
\begin{equation*}\label{arordam}
\begin{split}
&Nf = x - \left(g_1 - Nf_1\right) + \left(g_2 - Nf_2\right) \\
& \geq \alpha_*(x) + [p] - 3f_1 \\
& = N\alpha_*(f) + \alpha_*\left(\left(g_1 - Nf_1\right) - \left(g_2 -
  Nf_2\right)\right) + [p] -3f_1 \\
& \geq N\alpha_*(f) + [p] -3\left(f_1 + \alpha_*(f_2)\right) \succ
N\alpha_*(f),
\end{split}
\end{equation*}
and then from Lemma \ref{divisibility} that 
\begin{equation}\label{eqq12}
f \succ \alpha_*(f).
\end{equation} 
Since $N[e] \succ g_1 +
\alpha_*(g_2) \geq N\left(f_1 + \alpha_*(f_2)\right)$ we have also
that that $[e] \succ f_1 + \alpha_*(f_2)$. It follows then from Lemma
\ref{divisibility} that there are
projections $p,q \in B$ such that $[p] = f_1, [q] = f_2$ and $p +
\alpha(q) \leq e$.  Then $[\alpha(p)] +
[q]  \prec [p] + [\alpha(q)]$ by (\ref{eqq12}) and another application
of Lemma \ref{divisibility} implies that there is a partial
isometry $t \in M_2(B)$ such that
$$
t \left( \begin{matrix} \alpha(p) & 0 \\ 0 & q \end{matrix} \right)
t^* \leq \left( \begin{matrix} p + \alpha(q) & 0 \\ 0 & 0 \end{matrix}
\right) .
$$
Then
$$
t\left( \begin{matrix} \alpha(p) & 0 \\ 0 & 0 \end{matrix} \right) = \left( \begin{matrix} v & 0 \\ 0 & 0 \end{matrix} \right)
$$
and
$$
t\left( \begin{matrix} 0 & 1 \\ 1 & 0 \end{matrix} \right)\left(
  \begin{matrix} q & 0 \\ 0 &  0\end{matrix} \right) = \left(
  \begin{matrix} w & 0 \\ 0 & 0 \end{matrix} \right)
$$
where $w,v \in B$ are partial isometries such that $v^*v = \alpha(p)$,
$w^*w = q$ and $v\alpha(p)v^* + wqw^* < p + \alpha(q)$. Let
$u$ be the canonical unitary in the multiplier algebra of $B
\times_{\alpha} \mathbb Z$ which implements $\alpha$ on $B$. Set $s =
vup + wqu^* + (e-p-\alpha(q))$ and note that $s^*s = p + \alpha(q) +
(e-p-\alpha(q)) = e$ while $ss^* = v\alpha(p)v^* + wqw^* + (e - p -
\alpha(q)) < e$.
\end{proof}

Let $0 < \epsilon <\frac{1}{2} $ and define
$f^{\epsilon}_1,f^{\epsilon}_0 : [0,1] \to [0,1]$ such that
$$
f^{\epsilon}_1(t) = \begin{cases} 0, \ & 0  \leq t \leq
  \frac{\epsilon}{2} \\ 1, \ & t \in  [\epsilon,1] \\ \text{linear} , \ & \text{else} \end{cases}
$$
and
$$
f^{\epsilon}_0(t) = \max\{0, t - \epsilon\} .
$$

\begin{lemma}\label{lemma?} There is a $\delta > 0$ such that for all
  $\epsilon \in ]0,\frac{1}{2}[$ the following holds: When $b,b'$ are
  positive contractions in a $C^*$-algebra $B$ such that
 \begin{equation}\label{holdop}
\left\|f^{\epsilon}_1(b) -f^{\epsilon}_1(b')\right\| \leq \delta
\end{equation}
and $\overline{f^{\epsilon}_0(b')Bf^{\epsilon}_0(b')}$ contains a
  projection $p$, then
  $\overline{f^{\epsilon}_1(b)Bf^{\epsilon}_1(b)}$ contains a
  projection which is Murray-von Neumann equivalent to $p$.
\begin{proof} Let $\delta > 0$ be so small that $\overline{yBy}$
  contains a projection Murray-von Neumann equivalent to $q$ whenever
  $x,y,q$ are positive contractions in a $C^*$-algebra $B$ such that
  $q$ is a projection and $\left\|xyx -q\right\| \leq \delta$. This
  $\delta$ will work because $f_1^{\epsilon}(b')p = p$ and it follows
  therefore from (\ref{holdop}) that
  $\left\|pf^{\epsilon}_1(b)p -p\right\| \leq \delta$.
\end{proof}
\end{lemma}

\emph{Proof of Theorem \ref{semicrossed}}. The general setup for the proof is the following. Let $A_{\infty}$ be the inductive limit of the sequence
\begin{equation}\label{pr1} 
\begin{xymatrix}{
A \ar[r]^-{\beta} & A \ar[r]^{\beta} & A \ar[r]^{\beta} & \dots 
}\end{xymatrix}
\end{equation}
We can then define an automorphism $\alpha$ of $A_{\infty}$ such that
$\alpha \circ \rho_{\infty,n} = \rho_{\infty,n}
\circ \beta$,
where $\rho_{\infty,n} : A \to A_{\infty}$ is the canonical
$*$-homomorphism from the $n$'th level in the sequence (\ref{pr1})
into the inductive limit algebra. In this notation the inverse of
$\alpha$ is defined such that
$\alpha^{-1} \circ \rho_{\infty,n} =
\rho_{\infty,n+1}$. Let $e \in A_{\infty}$ be the projection $e =
\rho_{\infty,1}(1)$ which is a full projection of $A_{\infty}$ by
assumption i) and hence also a full projection of the crossed product
$A_{\infty} \times_{\alpha} \mathbb Z$. By a result of Stacey, \cite{St}, there
is an isomorphism $A \times_{\beta} \mathbb N \to e\left(A_{\infty} \times_{\alpha} \mathbb
  Z\right)e$ sending $a\in A$ to $\rho_{\infty,1}(a)$ and the
canonical isometry $v \in A \times_{\beta} \mathbb N$ to $eue$ where
$u$ is the canonical unitary in the multiplier algebra of $A_{\infty}
\times_{\alpha} \mathbb Z$. Note that $A_{\infty} \times_{\alpha}
\mathbb Z$ is stably isomorphic to $A \times_{\beta} \mathbb N$ and
hence simple by assumption. Thanks to condition i) the unit of
$\rho_{\infty,k}(A)$ is full in $\rho_{\infty,k+1}(A)$ so that the
sequence $B_k = \rho_{\infty,k}(A), k = 1,2,\dots $, will have
properties required in Lemma \ref{fullinf}. Furthermore, it follows
from condition ii) that there can not be any non-zero densely defined
lower semi-continuous $\alpha$-invariant trace on $A_{\infty}$;
because if there was it would have to be non-zero on some $B_k$ and it
would then give rise to a $\beta$-invariant trace state on $A$. In this
way it follows from Lemma \ref{fullinf} that every full projection of $A_{\infty}$ is
infinite in $A_{\infty}
\times_{\alpha} \mathbb Z$. In fact, the same argument shows that a
full projection in $M_k\left(A_{\infty}\right)$ is infinite in
$M_k\left(A_{\infty} \times_{\alpha} \mathbb Z\right)$ for any $k \in
\mathbb N$.

We make now the following

\begin{assert}\label{assert}
Let $h \in A_{\infty}\backslash \{0\}$ be a positive contraction. It follows that
  $\overline{h(A_{\infty} \times_{\alpha} \mathbb Z)h}$ contains an infinite projection.
\end{assert}

Assuming that Assertion \ref{assert} holds the proof of Theorem \ref{semicrossed} is completed as follows. Let $b \in \left(A_{\infty} \times_{\alpha} \mathbb Z\right) \backslash \{0\}$
be a positive contraction. Let $E : A_{\infty} \times_{\alpha} \mathbb
Z \to A_{\infty}$ be the canonical conditional expectation. Let
$\epsilon > 0$. As in the
proof of Lemma 2.4 of \cite{R2} we can find positive elements $ a,x \in
A_{\infty}$ such that $\left\|a\right\| \geq 1- \epsilon$, $\|x\| \leq 1$
and $\left\| \|E(b)\|^{-1} {xbx} - a\right\| \leq \epsilon$. The only
change we have to make to R\o rdams
argument is to replace the lemma of Kishomoto used by him with Lemma
7.1 of \cite{OP2}. Some backtracking through the work of Olesen and
Pedersen is needed to verify that Lemma 7.1 of \cite{OP2}
applies. What is needed is to show that the simplicity of $A_{\infty}
\times_{\alpha} \mathbb Z$ forces all the automorphisms $\alpha^n, n
\in \mathbb Z \backslash \{0\}$, to be properly outer since this is the
assumption in Lemma 7.1 of \cite{OP2}. This follows from the implication (i) $\Rightarrow$
(vi) of Theorem 10.4 in \cite{OP2} since the Connes spectrum
$\Gamma(\alpha)$ is the whole circle by Proposition 6.3 in \cite{OP1}.

 Having the element $a$, set $a' = f(a)$ where $f : [0,1] \to
[0,1]$ is a continuous function such that $f(t) = 1, t \in
[1-2\epsilon,1]$ and $|f(t)-t|\leq 2\epsilon$ for all $t \in
[0,1]$. Then $\left\|a' - a\right\| \leq 2 \epsilon$ and spectral
theory gives us a positive element $h \in A_{\infty}$ such that
$\left\|h\right\| = 1$ and $a'h = h$. It follows now from Assertion
\ref{assert} that $\overline{h\left(A_{\infty} \times_{\alpha} \mathbb
    Z\right)h}$ contains an infinite projection $p$. Since $\left\|
  \|E(b)\|^{-1}{xbx} - a'\right\| \leq 3\epsilon$ and $a'p = p$ we
find that $\left\|
 \|E(b)\|^{-1} {pxbxp}{} - p\right\| \leq 3\epsilon$. Thus, if only
$\epsilon$ is small enough
$\left\|E(b)\right\|^{-1}\sqrt{b}xpx\sqrt{b}$ will be close to a
projection in $\overline{b\left(A_{\infty} \times_{\alpha} \mathbb
    Z\right)b}$ which is Murray-von Neumann equivalent to $p$ and
hence infinite. This shows that $A_{\infty} \times_{\alpha} \mathbb Z$
is purely infinite, and the same is $A \times_{\beta} \mathbb N$ since
it is stably isomorphic to $A_{\infty} \times_{\alpha} \mathbb Z$,
cf. Proposition 5.5 of \cite{PS}. 

\smallskip 

It remains to prove Assertion \ref{assert}: Since $A_{\infty} =
\overline{\bigcup_n \rho_{\infty,n}(A)}$ an approximation argument
based on Lemma \ref{lemma?} shows that we may assume that $h \in
\rho_{\infty,n_0}(A)$ for some $n_0 \in \mathbb N$. Set $A'
= \rho_{\infty,n_0}(1)A_{\infty}\rho_{\infty,n_0}(1)$ and note that $uA'u^* \subseteq
A'$. 

%Furthermore, it follows from assumption i) about $\beta$ that
%$\beta^{k-n}(1)$ is a full projection in $A$ for all $k \geq n$. Since 
%$$
%A' = \overline{\bigcup_{k \geq n}  \rho_{\infty,n}(1)\rho_{\infty,k}(A)\rho_{\infty,n}(1)},
%$$
%it follows from Lemma \ref{AAA} that $A'$ is tracially almost
%divisible since we asssume that $A$ is.

Deviating slightly from the notation used so far, let $T(A_{\infty})$
denote the
  set of densely defined lower semi continuous traces $\omega$ on $A_{\infty}$ such that $\omega(e) =
  1$. This is a compact space in a topology described before Lemma 3
  of \cite{Th2} which is the same topology it gets through the
  identification of $T\left(A_{\infty}\right)$ with the tracial state
  space $T\left(eA_{\infty}e\right)$. Since $h \leq \rho_{\infty,n_0}(1)$ it follows that $\omega \mapsto \omega\left(u^kh{u^*}^k\right)$
is continuous on $T\left(A_{\infty}\right)$ for all $k$. We claim that there is an $m \in \mathbb N$ such that
\begin{equation}\label{B47}
\omega\left(h + uhu^* + u^2h{u^*}^2 + u^3h{u^*}^3 + \dots +
u^mh{u^*}^m\right) > 0
\end{equation}
for all $\omega \in T(A_{\infty})$. Indeed, if not there is for each $n \in
\mathbb N$ a trace $\omega_n \in T(A_{\infty})$ such that 
$$
\omega_n\left(h + uhu^* + u^2h{u^*}^2 + u^3h{u^*}^3 + \dots +
u^nh{u^*}^n\right) = 0.
$$
A condensation point of $\{\omega_n\}$ in $T(A_{\infty})$ will be a densely defined lower semi continuous trace $\omega$ such
that $\omega(u^nh{u^*}^n) = 0$ for all $n$. Then $\left\{ x \in A_{\infty} :
  \omega(u^nx^*x{u^*}^n) = 0 \ \forall n\right\}$ is a non-zero closed
two-sided ideal $I$ in $A_{\infty}$ such that $\alpha(I) \subseteq
I$. As in b)
from the proof of Lemma \ref{fullinf} this implies that $I = A_{\infty}$, which
is impossible since $e \notin I$. This proves the claim. 

Let $T(A')$ be the tracial state space of $A'$. Since $\beta(1)$ is full in $A$ it follows that $\rho_{\infty,n_0}(1)$ is a
full projection in $A_{\infty}$ so any $\omega \in T(A')$ is the
restriction to $A'$ of a densely defined lower semi-continuous trace
on $A_{\infty}$,
cf. Theorem 5.2.7 of \cite{Pe}. It follows therefore from (\ref{B47}) that
$$
\omega\left(h + uhu^* + u^2h{u^*}^2 + u^3h{u^*}^3 + \dots +
u^mh{u^*}^m\right) > 0
$$
for all $\omega \in T(A')$. By compactness of $T(A')$ there is a $\delta > 0$ such that\begin{equation*}\label{lowerest}
\omega\left(h + uhu^* + u^2h{u^*}^2 + u^3h{u^*}^3 + \dots +
u^mh{u^*}^m\right) \geq \delta
\end{equation*}
for all $\tau \in T(A')$. Since $A \simeq \rho_{\infty,n_0}(A)$ is tracially
almost divisible by Proposition \ref{tracAH} and since $\rho_{\infty,n_0}(A)$ is a
unital $C^*$-subalgebra of $A'$ which contains $h$ there
is a $\delta' > 0$ with the property that for any $\epsilon_1 > 0$ there are orthogonal
positive elements $h_1,h_2, \dots, h_{m+1}$ in $A'$ such that $h_1 +
h_2 + \dots + h_{m+1} \preceq h$ in $A'$ and $\tau(h_i) \geq \delta'\tau(h) -
\epsilon_1$ for all $i$ and all $\tau \in T(A')$. Let $\tau' \in  T\left(M_{m+1}(A')\right)$ - the tracial state space
of $M_{m+1}(A')$. Then $\tau' = \tau \otimes \tr$ for some $\tau \in
T(A')$, where $\tr$ is the trace state of $M_{m+1}$. 
It follows that
\begin{equation*}
\begin{split}
&\tau'  \left( \begin{matrix} h_1 & 0 & \hdots & 0\\ 0 & uh_2u^* & \hdots &  0\\
    \vdots & \vdots & \ddots & \vdots \\ 0 & 0 & \hdots & u^{m}h_{m+1}{u^*}^{m}
  \end{matrix}\right) =  \frac{1}{m+1} \sum_{j=1}^{m+1}
\tau\left(u^{j-1}h_j{u^*}^{j-1}\right)\\
&  \geq
\frac{1}{m+1} \sum_{j=1}^{m+1} \left(\delta'
{\tau\left(u^{j-1}h{u^*}^{j-1}\right)}
- \tau\left(\rho_{\infty,n}\left(\beta^{j-1}(1)\right)\right) \epsilon_1\right)  \\
& \geq \frac{\delta \delta'}{m+1} 
- \epsilon_1\left( \sum_{j=1}^{m+1} \tau\left(\rho_{\infty,n}\left(\beta^{j-1}(1)\right)\right)\right) .
\end{split}
\end{equation*}
Choose $\epsilon_1 > 0$
such that 
$$
\delta_1 = \frac{\delta \delta'}{m+1} 
- \epsilon_1 \sup_{\omega \in T(A')}\left( \sum_{j=1}^{m+1}
  \omega\left(\rho_{\infty,n}\left(\beta^{j-1}(1)\right)\right)\right) > 0.
$$ 
Set 
$$
H =  \left( \begin{matrix} h_1 & 0 & \hdots & 0\\ 0 & uh_2u^* & \hdots &  0\\
    \vdots & \vdots & \ddots & \vdots \\ 0 & 0 & \hdots & u^{m}h_{m+1}u^{m}
  \end{matrix}\right)
$$
and notice that $\tau(H) \geq \delta_1$ for all $\tau \in
T\left(M_{m+1}(A')\right)$. Let $\epsilon_0 \in ]0,\frac{1}{4}[$ be so
small that
\begin{equation}\label{laterref}
\tau\left(f^{\epsilon_0}_0(H)\right) \geq \frac{\delta_1}{2}
\end{equation}
for all $\tau \in T(M_{m+1}(A'))$.
It may or may not be the case that $A'$ is an AH-algebra with slow
dimension growth, but since $A$ has these properties and since
$$
A' = \overline{\bigcup_{k \geq n_0}
  \rho_{\infty,n_0}(1)\rho_{\infty,k}(A)\rho_{\infty,n_0}(1)} ,
$$
we can pick an
increasing sequence $ F_1 \subseteq F_2 \subseteq F_3 \subseteq \dots $
of finite subsets with dense union in $A'$ and write $A' =
\overline{\bigcup_l A_l}$ 
such that each $A_l$ is a homogeneous $C^*$-algebra with
$\rho_{\infty,n_0}(1) \in A_l$ and $F_l
\subseteq_{\frac{1}{l}} A_l$, meaning that every element of $F_l$ has
distance less than $\frac{1}{l}$ to an element of $A_l$, and such that $\lim_{l\to \infty} r(A_l)
= 0$. Let $\epsilon > 0$. We can then find $n_{\epsilon} \in \mathbb N$ and
for each $l \geq n_{\epsilon}$ a positive contraction $k_l \in M_{m+1}(A_l)$ such
that
\begin{equation}\label{laterref2}
\left\|f^{\epsilon_0}_1(H) - f^{\epsilon_0}_1\left(k_l\right)\right\| \leq
\epsilon
\end{equation}
and
$$
\left\|f^{\epsilon_0}_0(H) - f^{\epsilon_0}_0\left(k_l\right)\right\| \leq
\frac{\delta_1}{6}.
$$
In particular, it follows from the last condition and (\ref{laterref})
that there is an $n_{\epsilon}' \geq n_{\epsilon}$ such that 
\begin{equation}\label{tracelow}
\tau\left(f^{\epsilon_0}_0\left(k_l\right)\right) \geq \frac{\delta_1}{4}
\end{equation}
for all $\tau \in T(M_{m+1}(A_l))$ and all $l \geq
n'_{\epsilon}$. (This is proved
by contradiction. If there are arbitrary large $n_i$ for which
$T(M_{m+1}(A_{n_i}))$ contains an element $\tau_{n_i}$ with
$\tau_{n_i}\left(f^{\epsilon_0}_0\left(k_{n_i}\right)\right) <
\frac{\delta_1}{4}$, consider a state extension $\tilde{\tau_{n_i}}$ of
$\tau_{n_i}$ to $M_{m+1}(A')$. A weak* condensation point of
$\left\{\tilde{\tau_{n_i}}\right\}$ will be an element of
$T(M_{m+1}(A'))$ for which (\ref{laterref}) fails.)

Consider an $l \geq n'_{\epsilon}$. Let $d_{\tau}: M_{m+1}(A_l)^+ \to \mathbb R^+$ denote the dimension
function corresponding to $\tau \in T\left(M_{m+1}(A_l)\right)$,
i.e. $d_{\tau}(a) = \lim_{n \to \infty}
\tau\left(a^{\frac{1}{n}}\right)$. It follows then from
(\ref{tracelow})  
$$
d_{\tau}\left(f^{\epsilon_0}_0\left(k_{l}\right)\right) \geq
\frac{\delta_1}{4} .
$$
Since $\lim_{l \to \infty} r\left(M_{m+1}(A_l)\right) = 0$ it follows from well-known
properties of vector bundles, or from Theorem \ref{TOMS}, that for all large $l$ there is a projection $p_l \in M_{m+1}(A_l)$ with constant
rank 1 over the spectrum of $A_l$. Then $d_{\tau}\left(p_l\right) \leq
r\left(M_{m+1}(A_l)\right)$ and hence
$$
d_{\tau}\left(f^{\epsilon_0}_0\left(k_{l}\right)\right) \geq
\frac{\delta_1}{4} > 
\frac{r\left(M_{m+1}(A_l)\right)}{2} + d_{\tau}\left(p_l\right)
$$
for all $\tau \in T(M_{m+1}(A_l))$ when $l$ is large enough. Fix such
an $l$. Theorem
\ref{TOMS} gives us now a sequence $\{x_n\}$ in
$M_{m+1}(A_l)$ such that 
$$
\lim_{n \to \infty} x_n f^{\epsilon_0}_0\left(k_{l}\right)x_n^* =
p_l .
$$
Note that $p_l$ is a full projection in
$M_{m+1}(A_{\infty})$. As pointed out in the beginning of the proof
$p_l$ is then an infinite projection in $M_{m+1}\left(A_{\infty}
  \times_{\alpha} \mathbb Z\right)$. Note also that 
$\left\|x_n \sqrt{f^{\epsilon_0}_0\left(k_{l}\right)}\right\| \leq
2$
which combined with (\ref{laterref2}) implies that
$$
\left\|x_n
  \sqrt{f^{\epsilon_0}_0\left(k_{l}\right)}f^{\epsilon_0}_1(H)
  \sqrt{f^{\epsilon_0}_0\left(k_{l}\right)}x_n^* - x_n
  \sqrt{f^{\epsilon_0}_0\left(k_{l}\right)}f^{\epsilon_0}_1(k_{l})
  \sqrt{f^{\epsilon_0}_0\left(k_{l}\right)}x_n^*\right\| \leq 4
    \epsilon
$$
for all large $n$. Since $f^{\epsilon_0}_1f^{\epsilon_0}_0 =
f^{\epsilon_0}_0$ we conclude that
$$
\left\|x_n
  \sqrt{f^{\epsilon_0}_0\left(k_{l}\right)}f^{\epsilon_0}_1(H)
  \sqrt{f^{\epsilon_0}_0\left(k_{l}\right)}x_n^* - p_l\right\| \leq 5
    \epsilon
$$
for all large $n$. Let $X \in M_{m+1}(A_{\infty} \times_{\alpha} \mathbb Z)$ be the matrix
$$
X = \left( \begin{matrix} f^{\epsilon_0}_1\left(h_1\right) & 0  & \hdots
    & 0\\ uf^{\epsilon_0}_1\left(h_2\right)  & 0 & \hdots &  0\\
    \vdots & \vdots & \ddots & \vdots \\ u^{m}f^{\epsilon_0}_1\left(h_{m+1}\right) & 0 & \hdots & 0
  \end{matrix}\right) .
$$ 
Then $  XX^* = f^{\epsilon_0}_1(H)$ since the $h_i$'s are mutually
orthogonal and 
$$
X^*X = \diag \left( \sum_{j=1}^{m+1} f^{\epsilon_0}_1\left(h_j\right),
  0,0,\dots, 0\right),  
$$
i.e. $ f_1^{\epsilon_0}(H) \sim X^*X$ in the sense of \cite{KR} 
which implies that $f_1^{\epsilon_0}(H) \preceq X^*X$. Since $\sum_{j=1}^{m+1} f^{\epsilon_0}_1\left(h_j\right) =
f^{\epsilon_0}_1\left(\sum_{j=1}^{m+1} h_j\right)$ and $\sum_{j=1}^{m+1}
h_j \preceq h$ it follows from
Proposition 1.11 of \cite{W} that there is an element $h_0$ in the
hereditary $C^*$-subalgebra of $A'$ generated by $h$ such that
$\sum_{j=1}^{m+1} f^{\epsilon_0}_1\left(h_j\right) \preceq h_0$. Thus
$f^{\epsilon_0}_1(H) \preceq h'_0$, where $h'_0 = \diag (h_0,0,0,\dots, 0)$, i.e. there is a
sequence $\{z_n\}$ in $M_{m+1}(A_{\infty} \times_{\alpha} \mathbb Z)$ such that
$\lim_n z_nh'_0z_n^* = f_1^{\epsilon_0}(H)$. Then 
$$
\left\|x_n\sqrt{f^{\epsilon_0}_0\left(k_{l}\right)}z_{n'} h'_0z_{n'}^*\sqrt{f^{\epsilon_0}_0\left(k_{l}\right)}x_n^* - p_l \right\| \leq 6 \epsilon
$$
when $n'$ and $n$ are sufficiently large and it follows that 
$$
\sqrt{h'_0}z_{n'}^*\sqrt{f^{\epsilon_0}_0\left(k_{l}\right)}x_n^*x_n\sqrt{f^{\epsilon_0}_0\left(k_{l}\right)}z_{n'}^*\sqrt{h'_0}
$$ 
will be close to a
projection in $\overline{h'_0M_{m+1}\left(A_{\infty} \times_{\alpha} \mathbb Z\right)h'_0}$
which is equivalent to $p_l$. Since $p_l$ is infinite 
this gives us the desired projection, completing the proof of
Assertion \ref{assert} and hence also the proof of the theorem. \qed

\bibliographystyle{plain}
\bibliography{/Users/toke/Documents/Matematik/VigtigeArtikler}
\bibliography{/home/anne/tokemeie/Mac/Matematik/VigtigeArtikler}

\begin{thebibliography}{WWW} %antallet af W'er er bredeste ind


%\bibitem[AO]{AO} J. M. Aarts and L. Oversteegen, {\em Matchbox
 %   manifolds}, Continua (Cincinatti, OH, 1994), 3-14, Lecture Notes
 % in Pure and Appl. Math., 170, Dekker, New York, 1995. 

%\bibitem[ALNR]{ALNR} D. Adji, M. Laca, M. Nilsen and I. Raeburn, {\em
%    Crossed Products by Semigroups of Endomorphisms and the Toeplitz
%    Algebras of Ordered Groups}, Proc. Amer. Math. Soc. {\bf 122}
%  (1994), 1133-1141.


%\bibitem[AP]{AP} Adler, R. L. Adler and R. Palais, {\em Homeomorphic conjugacy of automorphisms on the torus},  Proc. Amer. Math. Soc.  {\bf 16} (1965), 1222--1225.


%\bibitem[AM]{AM} R. Adler and B. Marcus, {\em Topological entropy and equivalence of dynamical systems}, Mem. Amer. Math. Soc. {\bf 219} (1979).


\bibitem[An]{An} C. Anantharaman-Delaroche, {\em Purely infinite $C^*$-algebras arising from dynamical systems}, Bull. Soc. Math. France {\bf 125} (1997), 199--225.

%\bibitem[AR]{AR} C. Anantharaman-Delaroche and J. Renault, {\em Amenable Groupoids}, L'Enseignement Math\'ematique, Ge\'neve, 2000.


%\bibitem[Ao]{Ao} N. Aoki, {\em Expanding Maps of Solenoids},
%Mh. Math. {\bf 105 } (1988), 1--34.



%\bibitem[ABG]{ABG} R.J. Archbold, J.W. Bunce and K.D. Gregson, {\em Extensions of states of $C^*$-algebras, II}, Proc. Roy. Soc. Edinburgh, Sect.A {\bf 92} (1982), 113--122.


%\bibitem[Au]{Au} J. Auslander, {\em Transformation groups without
%    minimal sets}, Math. Systems Theory {\bf 2} (1967), 93-95.

%\bibitem[BEP]{BEP} T. Bates, S. Eilers and D. Pask, {\em Reducibility
%    of covers of AFT shifts}, arXiv:0812.1408

%\bibitem[AV]{AV} V.A. Azurmanian and A.M. Vershik, {\em Star algebras associated with endomorphisms}, in Operator Algebras and Group Representations, vol. {\bf I}, Ptiman, 1984, 17--27.


%\bibitem[BR1]{BR1} M. Baake and J. A. G. Roberts, {\em Trace maps as 3D reversible dynamical systems with an invariant}, J. Stat. Phys. {\bf 74} (1994), 829--888.


%\bibitem[BR2]{BR2} \bysame, {\em Reversing symmetry group of $Gl(2,\mathbb Z)$ and $PGl(2, \mathbb Z)$ matrices with connections to cat maps and trace maps}, J. Phys. A: Math. Gen. {\bf 30} (1997), 1549--1573.


%\bibitem[BEK]{BEK} O. Bratteli, D.E. Evans and A. Kishimoto, {\em
%    Crossed products of totally disconnected spaces by $\mathbb Z_2 *
%    \mathbb Z_2$}, Ergod. Th. \& Dynam. Sys. {\bf 13} (1993), 445-484.




%\bibitem[Ba]{Ba} C.J.K. Batty, {\em Simplexes of extensions of states of $C^*$-algebras}, Trans. Amer. Math. Soc. {\bf 272} (1982), 237--246.

%\bibitem[B1]{B1} S. Bhattacharya, {\em Orbit equivalence and topological conjugacy of affine actions on compact abelian groups}, Monatsh. Math.  {\bf 129}  (2000),  89--96.


%\bibitem[B2]{B2} \bysame, {\em Zero-entropy algebraic $\Bbb Z\sp d$-actions that do not exhibit rigidity},  Duke Math. J. {\bf 116} (2003), 471--476. 



%\bibitem[B3]{B3} \bysame, {\em Higher order mixing and rigidity of algebraic actions on compact abelian groups},  Israel J. Math. {\bf  137} (2003), 211--221. 

%\bibitem[BW]{BW} S. Bhattacharya and T. Ward, {\em Finite entropy characterizes topological rigidity},  Ergodic Theory and Dynamical Systems {\bf 25} (2005), 365--373.


%\bibitem[Bl]{Bl} B. Blackadar, {\em K-theory for Operator Algebras}, Springer Verlag, New York, 1986.

%\bibitem[BKR]{BKR} B. Blackadar, A. Kumjian and M. R\o rdam, {\em Approximately central matrix units and the structure of noncommative tori}, K-theory {\bf 6} (1992), 267--284.



\bibitem[BKR]{BKR} S. Boyd, N. Keswari and I. Raeburn, {\em Faithful
    Representations of Crossed Products by Endomorphisms},
  Proc. Amer. Math. Soc. {\bf 118} (1993), 427-436.

%\bibitem[B]{B} M. Boyle, {\em Lower entropy factors of sofic systems}, Ergod. Th. \& Dynam. Sys. {\bf 4} (1984), 541--557. 


%\bibitem[BFF]{BFF} M. Boyle, D. Fiebig and U. Fiebig, {\em A dimension group for local homeomorphisms and endomorphisms of onesided shifts of finite type},  J. Reine Angew. Math.  {\bf 487} (1997), 27--59.



%\bibitem[B-M]{B-M} A. Bertrand-Mathis, {\em D\'eveloppement en base $\theta$, r\'epartition modulo un de la suite $(x\theta^n)_{n \geq 0}$, langages codes et $\theta$-shift}, Bull. Soc. math. France {\bf 114}(1986), 271--323.


%\bibitem[B-H]{B-H} Y. Bar-Hilel, M. Peres and E. Shamir, {\em On formal properties of simple phrase structure grammar} Z. Phonetik. Sprachwiss. Komm. {\bf 14} (1961), 143--172.


 
%\bibitem[BS]{BS} S. Bhattacharya and K. Schmidt, {\em Homoclinic points and isomorphism rigidity of algebraic $\Bbb Z\sp d$-actions on zero-dimensional compact abelian groups},  Israel J. Math. {\bf 137} (2003), 189--209.

%\bibitem[Bl]{Bl} F. Blanchard, {\em $\beta$-expansions and symbolic dynamics}, Theoretical Computer Science {\bf 65} (1989) 131--141.


%\bibitem[BH]{BH} F. Blanchard and G. Hansel, {\em Syst\`emes code\'es}, Theoret. Comp. Sci. {\bf 44} (1986), 17--49.


%\bibitem[BK]{BK} M. Boyle and W. Krieger, {\em Almost Markov and shift equivalent sofic systems}, Proceedings of the Maryland Special Year in Dynamics 1986--87, Springer-Verlag, LNM {\bf 1342} (1988), 33--93.

%\bibitem[BBG]{BBG} M. Boyle, J. Buzzi and R. Gomez, {\em Almost isomorphisms for countable state Markov shifts}, Preprint, 2004.


%\bibitem[BKM]{BKM} M. Boyle, B. Kitchens and B. Marcus, {\em A note on minimal covers for sofic systems}, Proc. Amer. Math. Soc. {\bf 95} (1985), 403--411.

%\bibitem[BFF]{BFF} M. Boyle, D. Fiebig and U. Fiebig, {\em Residual entropy, conditional entropy and subshift covers}, Preprint (2000).

%\bibitem[Bo]{Bo} R. Bowen, {\em Markov partitions for axiom A diffeomorphisms}, Amer. J. Math. {\bf 92} (1970), 725--747.

%\bibitem[B1]{B1} R. Bowen, {\em Topological Entropy and Axiom A}, Proceedings of Symposia in Pure Mathematics {\bf 14} (1970), 23--41.


%\bibitem[B2]{B2} \bysame, {\em Entropy for group endomorphisms and homogeneous spaces}, Trans. Amer. Math. Soc. {\bf 153} (1971), 401--414.


%\bibitem[Bu]{Bu} \bysame, {\em Subshifts of Quasi-Finite Type}, Invent. Math., to appear.


%\bibitem[BJKR]{BJKR} O. Bratteli, P. E. T. Jorgensen, K.H. Kim and
%  F. Roush {\em Computation of isomorphism invariants for stationary
%    dimension groups},  Ergod. Th. \& Dynam. Sys. {\bf 22} (2002), 99--127.

%\bibitem[Bre]{Bre} B. Brenken, {\em The local product structure of expansive automorphisms of solenoids and their associated $C^*$-algebras}, Can. J. Math. {\bf 48} (1996), 692--709.


%\bibitem[Br]{Br} L. Brown, {\em Stable isomorphism of hereditary subalgebras of $C^*$-algebras}, Pacific J. Math. {\bf 71} (1977), 335--348.


%\bibitem[BP]{BP} L. Brown and G.K. Pedersen, {\em $C^*$-algebras of real
%    rank zero}, J. Funct. Anal. {\bf 99} (1991), 132--149.


%\bibitem[BGR]{BGR} L. Brown, P. Green, M. Rieffel, {\em Stable isomorphism and strong Morita equivalence of $C^*$-algebras}, Pacific J. Math. {\bf 71} (1977), 349--363.

%\bibitem[C]{C} T.M. Carlsen, {\em Cuntz-Pimsner $C^*$-algebras
%    associated with subshifts}, Internat. J. Math. {\bf 19} (2008), 47-70.
\bibitem[CT]{CT} T.M. Carlsen and K. Thomsen, {\em The structure of
    the $C^*$-algebra of a locally injective surjection}, in preparation.

%\bibitem[CS]{CS} T.M. Carlsen and S. Silvestrov, {\em $C^*$-crossed
%    products and shift spaces}, Expo. Math. {\bf 25} (2007), 275-307.

%\bibitem[CM]{CM} T.M. Carlsen and K. Matsumoto, {\em Some remarks on
%    the $C^*$-algebras associated with subshifts}, Math. Scand. {\bf
%    95} (2004), 145-160.


 %\bibitem[CE]{CE} M. D. Choi and G. Elliott, {\em Density of the
 %   selfadjoint elements with finite spectrum in an irrational
  %  rotation $C\sp *$-algebra}  Math. Scand. {\em  67} (1990), 73--86.

%\bibitem[CF]{CF} A. Clark and R. Fokkink, {\em On a Homoclinic Group that is not Isomorphic to the Character Group}, Qualitative Theory of Dynamical Systems {\bf 5} (2004), 361--365.


%\bibitem[C]{C} D. L. Cohn, {\em Measure Theory}, Birkhäuser, Boston, 1980.

%\bibitem[Co]{Co} A. Connes, {\em A survey of foliations and operator algebras}, Proc. Symp. Pure Math. {\bf 38(1)}, (1982), 521--628.

%\bibitem[CV]{CV} E.M. Coven, W.L. Reddy, {\em Positively expansive maps of compact manifolds},  Global theory of dynamical systems (Proc. Internat. Conf., Northwestern Univ., Evanston, Ill., 1979), pp. 96--110, Lecture Notes in Math., 819, Springer, Berlin, 1980.



%\bibitem[CK]{CK} E. M. Coven and M. Keane {\em Every compact metric psace that supports a positively expansive homeomorphism is finite},  IMS Lecture Notes--Monograph Series Vol. 48 (2006), 304--305.




%\bibitem[CP]{CP} E. Coven and M. Paul, {\em Finite procedures for sofic systems}, Monats. Math. {\bf 83} (1977), 265--278.

%\bibitem[CuPe]{CuPe} J. Cuntz and G.K. Pedersen, {\em Equivalence and
%    trace in $C^*$-algebas}, JFA, 1979.

%\bibitem[CK]{CK} J. Cuntz and W. Krieger, {\em A class of
%    $C^*$-algebras and topological Markov chains}, Invent. Math. {\bf
%    56} (1980), 25--268. 

\bibitem[De]{De} V. Deaconu, {\em Groupoids associated with endomorphisms}, Trans. Amer. Math. Soc. {\bf 347} (1995), 1779-1786.


%\bibitem[DS]{DS} V. Deaconu and F. Schultz, {\em $C^*$-algebras
%    associated with interval maps}, Trans. Amer. Math. Soce. {\bf 359}
%  (2007), 1889-1924.



%\bibitem[Di]{Di} J. Dixmier, {\em Von Neumann algebras}, North-Holland Publishing company, Amsterdam, New York, Oxford, 1981.

%\bibitem[Di]{Di} J. Dixmier, {\em $C^*$-algebras}, North-Holland Publishing company, Amsterdam, New York, Oxford, 1977.



%\bibitem[DT1]{DT1} S. Dorofeev and K. Thomsen, {\em Factors and subfactors arising from inductive limits of interval algebras}, Ergod. Th. \& Dynam. Sys. {\bf 19} (1999), 363--381.
 

%\bibitem[DT2]{DT2} \bysame , {\em Factors from ergodic theory and group-subgroup subfactors}, J. Ramanujan Math. Soc. {\bf 12} (1997), 239--262.


%\bibitem[DHS]{DHS} F. Durand, B. Host and C. Skau, {\em Substitutional
 %   dynamical systems, Bratteli diagrams and dimension groups},
 % Ergod. Th. \& Dynam. Syst. {\bf 19} (1999), 953--993.


\bibitem[DR]{DR} K. Dykema and M. R\o rdam, {\em Purely infinite,
    simple $C^*$-algebras arising from free product constructions},
  Canad. J. Math. {\bf 50} (1998), 323-341.  

%\bibitem[Ef]{Ef} E. Effros, {\em Dimensions and $C^*$-algebras}, CBMS
%  Regional Conf. Ser. in Math., no. 46, Amer. Math. Soc., Providence, 1981.


%\bibitem[EO]{EO} G. Elliott and D. Olesen, {\em A simple proof of the
%    Dauns-Hofmann theorem}, Math. Scand. {\bf 34} (1974), 231-234.



%\bibitem[Ell]{Ell} G. Elliott, {\em On the K-theory of the $C^*$-algebra generated by a projective representation of a torsion-free discrete abelian group}, pages 157--184 in: {\em Operator Algebras and Group Representations, Vol. I (neptun, 1980)}, Monogr. Stud. Math. {\bf 17}, Pitman, Boston MA, 1984.


%\bibitem[Ell1]{Ell1}  G. Elliott, {\em A classification of certain simple $C^*$-algebras, II}, J. Ramanujan Math. Soc. {\bf 12} (1997), 97--134.

%\bibitem[Ell2]{Ell2} \bysame, {\em On the classification of $C\sp
%    *$-algebras of real rank zero},  J. Reine Angew. Math. {\bf 443} (1993),
%  179--219.



%\bibitem[Ell2]{Ell2} \bysame, {\em  An invariant for simple $C\sp *$-algebras},  Canadian Mathematical Society. 1945--1995, Vol. 3,  61--90, Canadian Math. Soc., Ottawa, ON, 1996. 


%\bibitem[EE]{EE} G. Elliott, D. Evans, {\em The structure of the irrational rotation $C^*$-algebra}, Ann. of Math. (2) {\bf 138} (1993), 477--501.


%\bibitem[EG]{EG} G. Elliott and G. Gong, {\em On the classification of
%    $C^*$-algebras of real rank zero, II}, Ann. Math. {\bf 144}
%  (1996), 497--610.

%\bibitem[EGL]{EGL} G. Elliott, G. Gong and L. Li, {\em On the
%    classification of simple inductive limit $C\sp
%    *$-algebras. II. The isomorphism theorem},
%  Invent. Math. {\bf 168} (2007), no. 2, 249-320.

%\bibitem[En]{En} R. Engelking, {\em Dimension Theory}, North-Holland, Amsterdam, Oxford, New York (1978)

%\bibitem[EW]{EW} M. Einsiedler and T. Ward, {\em Thomas Isomorphism rigidity in entropy rank two},  Israel J. Math. {\bf  147} (2005), 269--284.


%\bibitem[ER]{ER} R. Exel and J. Renault, {\em Semigroups of local
%    homeomorphisms and interaction groups}, Ergod. Th. \&
%  Dynam. Sys. {\bf 27} (2007), 1737-1771.

%\bibitem[EV]{EV} R. Exel and A. Vershik, {\em $C^*$-algebras of
%    Irreversible Dynamical Systems}, Canad. J. Math. {\bf 58} (2006), 39-63.


%\bibitem[Fo]{Fo} R. Fokkink, {\em The structure of trajectories},
%  Ph.D. Thesis, Technische Universiteit te Delft, 1991.

%\bibitem[F]{F} D. Fiebig, {\em Factor maps, entropy and fiber cardinality for Markov shifts}, Rocky Mt. J. Math. {\bf 31} (2001), 955--986.

%\bibitem[FF]{FF} D. Fiebig and U. Fiebig, {\em Topological boundaries for countable state Markov shifts}, Proc. London Math. Soc., Ser. {\bf III 70} (1995)  , 625--643.

%\bibitem[FF]{FF} D. Fiebig and U. Fiebig, {\em Covers for coded systems}, in {\em Symbolic Dynamics and Its Applications}, Contemporary Mathematics {\bf 135} (ed. P. Walters), Providence, 1992, pp. 139--180.


%\bibitem[FFJ]{FFJ} D. Fiebig, U. Fiebig, N. Jonoska, {\em Multiplicities of covers for sofic shifts}, Theor. Comp. Science {\bf 262} (2001), 349--375.


%\bibitem[FR]{FR} J. Franks and D. Richeson, {\em Shift equivalence and the Conley index}, Trans. Amer. Math. Soc. {\bf 352} (2000), 3305--3322.

%\bibitem[Fu]{Fu} L. Fuchs, {\em Abelian groups}, Pergamon Press, Oxford,
%  London, New York, Paris, 1960.


%\bibitem[GPS1]{GPS1} T. Giordano, I. Putnam and C. Skau, {\em
%    Topological orbit equivalence and $C^*$-crossed products}, Journal
%  f\"ur die reine und angewandte Mathematik {\bf 469} (1995), 51--111.


%\bibitem[GPS2]{GPS2} \bysame, {\em Affable equivalence relations and orbit structure of Cantor dynamical systems},  Ergodic Theory Dynam. Systems  {\bf 24}  (2004),  no. 2, 441--475.



%\bibitem[Go]{Go} G. Gong, {\em On the classification of Simple Inductive Limit $C^*$-algebras, I: The Reduction Theorem}, Documenta Math. {\bf 7} (2002), 255--461.

%\bibitem[G]{G} K. R. Goodearl, {\em Partially ordered abelian groups with interpolation}, Mathematical Surveys and Monographs, 20. American Mathematical Society, Providence, RI, 1986. 



%\bibitem[GD]{GD} A. Granas and J. Dugundji, {\em Fixed Point Theory},
%  Springer Verlag, New York, 2003. 

%\bibitem[GHT]{GHT} E. Guentner, N. Higson and J. Trout, {\em Equivariant $E$-theory for $C\sp *$-algebras},  Mem. Amer. Math. Soc.  148  (2000),  no. 703.


 
%\bibitem[Gu1]{Gu1} B.M. Gurevi\v{c}, {\em Topological entropy of enumerable Markov chains} (Russian), Dokl. Akad. Nauk SSSR {\bf 187} (1969), 715--718; English translation: Sovjet Math. Dokl. {\bf 10(4)} (1969), 911--915.


%\bibitem[Gu2]{Gu2} \bysame, {\em Shift entropy and Markov measures in the path space of a denumerable graph} (Russian), Dokl. Akad. Nauk SSSR {\bf 192} (1970), 963--965; English translation: Sovjet Math. Dokl. {\bf 11(3)} (1970), 744--747.

%\bibitem[GS]{GS} B.M. Gurevi\v{c} and S.V. Savchenko, {\em Thermodynamic formalism for countable symbolic Markov chains} (Russian), Uspekhi Mat. Nauk {\bf 53(2)}, 3--106; English translation: Russian Math. Surveys {\bf 53(2)}, 245--344.




%\bibitem[GZ]{GZ} B.M. Gurevi\v{c} and A.S. Zargaryan, {\em Existence conditions of a maximal measure for a countable symbolic Markov chain} (Russian), Vestnik Moskov. Univ. Ser. I Mat. Mekh., {\bf 43(5)} (1988), 14--18; English translation: Moscow Univ. Math. Bull., 1988, 18--23.




%\bibitem[HPS]{HPS} R.H. Herman, I.F. Putnam and C.F. Skau, {\em
%    Ordered Bratteli diagrams, dimension groups and topological
%    dynamics}, Internat. J. Math. {\bf 3} (1992), 827--864.


%\bibitem[HI]{HI} T. Harju and L. Ilie, {\em Languages obtained from infinite words},  RAIRO Inform. Théor. Appl.  {\bf 31}  (1997), 445--455.

%\bibitem[HK]{HK} B. Hasselblatt and A. Katok, {\em Introduction to the Modern Theory of Dynamical Systems}, Encyclopedia of Math. and its Appl. {\bf 54}, Cambridge Univ. Press (1995).

%\bibitem[HeRo]{HeRo} E. Hewitt and K. Ross, {\em Abstract Harmonic Analysis I}, Springer Verlag, Berlin, G\"ottingen, Heidelberg, 1963.

%\bibitem[HSZ]{HSZ} M.W. Hirsch, H.L. Smith and Xiao-Qiang Zhao, {\em    Chain transitivity, Attractivity, and Strong Repellors for Semidynamical Systems}, J. Dyn. Diff. Eq. {\bf 13} (2001), 107--131.




%\bibitem[HR]{HR} J. v. B. Hjelmborg, M. R\o rdam, {\em On stability of $C^*$-algebras}, J. Funct. Anal. {\bf 155} (1998), 152--170.




 
%\bibitem[HU]{HU} J. E. Hopcroft and J. D. Ullman, {\em Introduction to Automata Theory, Languages and Computation}, Addison-Wesley Publishing Company, 1979.


%\bibitem[JS]{JS} X. Jiang and H. Su, {\em On a simple unital projectionless $C^*$-algebra}, Amer. J. Math. {\bf 121} (1999), 359--413.


%\bibitem[J]{J} V. Jones, {\em Index for subfactors}, Invent. Math. {\bf 72} (1983), 1--25.

%\bibitem[J]{J} N. Jonoska, {\em Sofic shifts with synchronizing presentations}, Theor. Comp. Science {\bf 158} (1996), 81--115.
 
%\bibitem[KPS]{KPS} J. Kaminker, I. Putnam and J. Spielberg, {\em Operator algebras and hyperbolic dynamics}, Operator Algebras and Quantum Field Theory, S. Doplicher, R. Longo, J.E.Roberts and L. Zsido, Eds., International Press, 1997.



%\bibitem[K]{K} T. Katsura, {\em A class of $C^*$-algebras generalizing
%    both graph algebras and homeomorphism algebras. III Ideal
%    structures}, Ergodic Th. Dyn. Syst. {\bf 26} (2006), 1805-1854.

%\bibitem[Ka1]{Ka1} G.G. Kasparov, J Op. Th.

%\bibitem[Ka2]{Ka2} \bysame , {\em Equivariant KK-theory and the Novikov conjecture}, Invent. Math. {\bf 91} (1988), 513--572.


 \bibitem[KR]{KR} E. Kirchberg and M. R\o rdam, {\em Non-simple purely
     infinite $C^*$-algebras}, Amer. J. Math. {\bf 122} (2000), 637-666.

%\bibitem[K1]{K1} B.P. Kitchens, {\em The dynamics of group automorphisms}, Annals of Numerical Mathematics {\bf 4} (1997), 369--391.

%\bibitem[K]{K}  B.P. Kitchens,  {Symbolic Dynamics}, Springer Verlag, Berlin, Heidelberg, 1998.

%\bibitem[KS]{KS}  B. Kitchens and K. Schmidt, {\em  Automorphisms of compact groups}, Ergodic Theory Dynam. Systems  {\bf 9}(1989), 691--735. 


%\bibitem[KS2]{KS2} \bysame {\em  Isomorphism rigidity of irreducible algebraic $Z\sp d$-actions},  Invent. Math.  {\bf 142}  (2000),  no. 3, 559--577.


%\bibitem[Ki]{Ki} B.P. Kitchen, {\em Symbolic Dynamics}, Springer Verlag, Berlin, Heidelberg, 1998.

%\bibitem[K-JT]{K-JT} K. Knudsen-Jensen and K. Thomsen, {\em Elements of KK-theory}, Birkhauser, Bosten, 1990.


%\bibitem[Kr1]{Kr1} W. Krieger, {\em  On a dimension for a class of homeomorphism groups}, Math. Ann. {\bf 252} (1979/80), 87--95. 


%\bibitem[Kr2]{Kr2} \bysame, {\em On dimension functions and topological Markov chains}, Invent. Math. {\bf 56} (1980), 239--250. 



%\bibitem[K]{K} W. Krieger, {\em On the Subsystems of Topological Markov Chains}, Ergod. Th. \& Dynam. Sys. {\bf 2} (1982), 195--202. 

%\bibitem[Kr]{Kr} W. Krieger, {\em On the uniqueness of the equilibrium state}, Math. Systems Theory {\bf 8} (1974), 97--104.
 
%\bibitem[Ku]{Ku} A. Kumjian, {\em On $C\sp *$-diagonals},  Canad. J. Math. {\bf 38}(1986), 969--1008.

%\bibitem[KuR]{KuR} A. Kumjian and J. Renault, {\em KMS-states on $C^*$-algebras associated to expansive maps}, Proc. Amer. Math. Soc. {\bf 134} (2006), 2067--2078.

 
%\bibitem[KPRR]{KPRR} A. Kumjian, D. Pask, I. Raeburn, J. Renault, {\em Graphs, groupoids, and Cuntz-Krieger algebras},  J. Funct. Anal.  {\bf 144}  (1997),  no. 2, 505--541.


%\bibitem[L]{L} P.-F. Lam, {\em On expansive transformation groups}, Trans. Amer. Math. Soc. {\bf 150} (1970), 131--138.


%\bibitem[La]{La} W. Lawton, {\em The structure of compact connected groups which admit an expansive automorphism}, Recent Advances in Topological Dynamics, LNS 318, Springer Verlag, 1973, 182--196.

%\bibitem[L1]{L1} H. Lin,  {\em Almost commuting selfadjoint matrices and applications},  Operator algebras and their applications (Waterloo, ON, 1994/1995),  193--233, Fields Inst. Commun., 13, Amer. Math. Soc., Providence, RI, 1997.


%\bibitem[Lin1]{Lin1} H. Lin, {\em Tracially AF $C^*$-algebras},
%  Trans. Amer. Math. Soc. {\bf 353} (2001), 693-722.

%\bibitem[Lin2]{Lin2}\bysame, {\em The tracial topological rank of $C^*$-algbras},
%Proc. London Math. Soc. {\bf 83} (2001), 199-234.


%\bibitem[Lin3]{Lin3} \bysame, {\em Traces and simple $C^*$-algebras with tracial topological rank zero}, J. Reine angew. Math. {\bf 568} (2004), 99--137.

%\bibitem[Lin4]{Lin4} \bysame, {\em Classification of simple
%    $C^*$-algebras of tracial rank zero}, Duke Math. J. {\bf 125}
%  (2004), 91-119.


%\bibitem[LM]{LM} D. Lind and B. Marcus, {\em An Introduction to
%    Symbolic Dynamics and Coding}, Cambridge University Press, 1995. 

%\bibitem[LS]{LS} D. Lind and K. Schmidt, {\em Homoclinic points of algebraic $\mathbb Z^d$-ations}, J. Amer. Math, Soc. {\bf 12} (1999), 953--980. 


%\bibitem[M]{M} R. Man\'e, {\em Expansive homeomorphisms and topological dimension},  Trans. Amer. Math. Soc. {\bf 252}(1979), 313--319.
 
%\bibitem[MRW]{MRW} P. Muhly, J. Renault and D. Williams, {\em Equivalence and isomorphism for groupoid $C^*$-algebras}, J. Oper. Th. {\bf 17} (1987), 3--22.

%\bibitem[M]{M} B. Marcus, {\em Sofic systems and encoding data}, IEEE Trans. Inform. Theory {\bf 31} (1985), 366--377.

%\bibitem[MT]{MT} J. Milnor and W. Thurston, {\em On iterated maps of the interval}, LNM 1342, pp. 465--563, Springer Verlag, Berlin, 1988.  


%\bibitem[Mu]{Mu} G.J. Murphy, {\em Simplicity of crossed products by
%      endomorphisms}, Integr. equ. oper. theory {\bf 42} (2002), 90-98. 

%\bibitem[N]{N} T. Natsume, {\em  On $K_*\left(C^*\left(SL_2(\mathbb
%        Z)\right)\right)$}, J. Operator Theory {\bf 13} (1985),
%  103-118.



%\bibitem[N]{N} M. Nasu, {\em An invariant for bounded-to-one factor maps between transitive sofic subshifts}, Ergod. Th. \& Dynam. Sys. {\bf 3} (1985), 89--105.

%\bibitem[N2]{N2} \bysame, {\em Topological conjugacy for sofic systems and extensions of automorphisms of finite subsystems of topological Markov shiftsA}, Proceedings of the Maryland Special Year in Dynamics 1986--87, Springer-Verlag, LNM {\bf 1342} (1988), 564--607.




%\bibitem[NP]{NP} D. Newton and W. Parry, {\em On a factor automorphism of a normal dynamical system}, Ann. Math. Statist. {\bf 37} (1966), 1528--1533.



%\bibitem[MP]{MP} S. Marcus and G. Paun, {\em Infinite (almost periodic) words, formal languages and dynamical systems}, Bulletin of the EATCS {\bf 54} (1994), 224--231.


%\bibitem[Ma]{Ma} K. Matsumoto, {\em $C^*$-algebras associated
%    with presentations of subshifts}, Doc. Math. {\bf 7} (2002), 1-30.


%\bibitem[Ma1]{Ma1} K. Matsumoto, {\em On $C^*$-algebras associated
%    with subshifts}, Internat. J. Math. {\bf 8} (1997), 357-374

%\bibitem[Ma2]{Ma2} \bysame , {\em $K$-theory for $C^*$-algebras
%    associated with subshifts}, Math. Scand. {\bf 82} (1998), 237-255.

%\bibitem[Ma3]{Ma3} \bysame, {\em Relation among generators of
%    $C^*$-algebras associated with subshifts}, Internat. J. Math. {\bf
%    10} (1999), 385-405.

%\bibitem[Ma4]{Ma4} \bysame, {\em On automorphisms of $C^*$-algebras
%    associated with subshifts}, J. Operator Theory {\bf 44} (2000),
%  81-112.

%\bibitem[Ma5]{Ma5} \bysame, {\em Stabilized $C^*$-algebras constructed
%    from symbolic dynamical systems}, Ergodic Theory Dynam. Systems
%  {\bf 20} (2000), 91-112.


%\bibitem[MW]{MW} P. Muhly and D. Williams, {\em Continuous trace groupoid $C^*$-algebras}, Math. Scand. {\bf 66}(1990), 239--250.

%\bibitem[R]{R} A. R\'enyi, {\em Representations for real numbers and their ergodic properties}, Acta. Math. Acad. Sci. Hungar. {\bf 8} (1957), 401--414.

%\bibitem[OP]{OP} H. Osaka and N. C. Phillips, {\em Crossed products by
%    finite groups with the Rokhlin property}, Preprint (2007),
%    arXiv:math/0704.3651v [math.OA]

\bibitem[OP1]{OP1} D. Olesen and G.K. Pedersen, {\em Applications of
    the Connes spectrum to $C^*$-dynamical systems}, J. Func. Analysis {\bf 30} (1978), 179-197.

\bibitem[OP2]{OP2} \bysame, {\em Applications of
    the Connes spectrum to $C^*$-dynamical systems, III},
  J. Func. Analysis {\bf 45} (1982), 357-390.

%\bibitem[P]{P} W.L. Paschke, {\em The Crossed Product of a
%    $C^*$-algebra by an Endomorphism}, Proc. Amer. Math. Soc. {\bf
%    80} (1980), 113-118.

%\bibitem[Pe]{Pe} G. K. Pedersen, {\em $C^*$-algebras and their automorphism groups}, Academic Press, London, New York, San Francisco, 1979.


%\bibitem[PR]{PR} J. Packer and I. Raeburn, {\em Twisted crossed products of $C^*$-algebras}, Math. Proc. Camb. Phil. Soc. {\bf 106} (1989), 293--311.

%\bibitem[Pa]{Pa} A.L.T. Paterson, {\em Groupoids, Inverse Semigroups, and their Operator Algebras}, Progress in Mathematics 170, Birkh\"auser, Berlin, 1999.


%\bibitem[P1]{P1} W.Parry, {\em On the $\beta$-expansion of real numbers}, Acta. Math. Acad. Sci. Hungar. {\bf 11} (1960), 4014--16.

%\bibitem[P2]{P2} \bysame, {\em Representations for real numbers}, Acta. Math. Acad. Sci. Hungar. {\bf 15} (1964), 95--105.


%\bibitem[Pa]{Pa} W. Parry, {\em Symbolic dynamics and transformations of the unit interval}, Trans. Amer. Math. Soc. {\bf 122} (1966), 368--378. 



%\bibitem[P]{P} W. Parry, {\em A finitary classification of topological Markov chains and sofic systems}, Bull. London Math. Soc. {\bf 9} (1977), 86--92.


%\bibitem[Pe]{Pe} K. Petersen, {\em Chains, entropy, coding}, Ergod. Th. \& Dynam. Sys. {\bf 6} (1986), 415--448.


\bibitem[Pe]{Pe} G. K. Pedersen, {\em $C^*$-algebras and Their Automorphism Groups}, London Mathematical Society Monographs, Vol. 14 (London: Academic Press, 1979). 


%\bibitem[Ph]{Ph} N. C. Phillips, {\em Cancellation and stable rank for
%    direct limits of recursive subhomogeneous $C^*$-algebras},
%  Trans. Amer. Math. Soc. {\bf 359} (2007), 4625-4652.



%\bibitem[Ph1]{Ph1} N. C. Phillips, {\em Crossed Products of the Cantor Set by Free Minimal Actions of $\mathbb Z^d$}, Comm. Math. Phys. {\bf 256} (2005), 1--42.

%\bibitem[Ph2]{Ph2} \bysame, {\em Every simple higher dimensional noncommutative torus is an AT algebra}, Preprint, September, 2006.

%\bibitem[Ph3]{Ph3} \bysame, {\em Real rank and property (SP) for
%    direct limits of recursive subhomogeneous algebras},
%  Trans. Amer. Math. Soc., to appear.

%\bibitem[Ph4]{Ph4} \bysame, {\em The tracial Rokhlin property
%    for actions of finite groups on $C^*$-algebras}, Preprint (2006), arXiv:math/0609782v1 [math.OA]


%\bibitem[Ph5]{Ph5} \bysame, {\em Finite cyclic group actions with the
%    tracial Rokhlin property}, Preprint (2006), arXiv:math/0609785v1 [math.OA]

%\bibitem[PV]{PV} M. Pimsner and D. Voiculescu, {\em Exact sequences for $K$-groups and $\Ext$-groups of certain cross-products of $C^*$-algebras}, J. Operator Theory {\bf 4} (1980), 93--118.



%\bibitem[Pu]{Pu} I. Putnam, {\em $C^*$-algebras from Smale spaces}, Can. J. Math. {\bf 48} (1996), 175--195.


%\bibitem[Pu2]{Pu2} \bysame, {\em Hyperbolic systems and generalized Cuntz-Krieger algebras}, Lecture notes from a Summer School in Operator Algebras, Odense, Denmark, 1996.

%\bibitem[Pu3]{Pu3} \bysame, {\em Functoriality of the $C^*$-algebras associated with hyperbolic dynamical systems}, J. London Math. Soc. (2) {\bf 62} (2000), 873--884.




%\bibitem[Pu4]{Pu4} \bysame, {\em Lifting factor maps to resolving maps}, Israel J. Math. {\bf 146}(2005), 253--280.




\bibitem[PS]{PS} I. Putnam and J. Spielberg, {\em The structure of $C^*$-algebras associated with hyperbolic dynamical systems}, J. Func. Anal. {\bf 163} (1999), 279--299.

%\bibitem[RW]{RW} I. Raeburn and D. P. Williams, {\em Morita
%    Equivalence and Continuous-Trace $C^*$-algebras}, American
%  Mathematical Society, 1998.


%\bibitem[Rd]{Rd} W. L. Reddy {\em Expanding maps on compact metric spaces}, Topology and its Applications {\bf 13} (1982), 327--334.


%\bibitem[Re]{Re} J. Renault, {\em A Groupoid Approach to $C^*$-algebras},  LNM 793, Springer Verlag, Berlin, Heidelberg, New York, 1980.  

%\bibitem[Re2]{Re2} \bysame, {\em The ideal structure of groupoid
%    crosssed product $C^*$-algebras}, J. Operator Theory {\bf 25}
%  (1991), 3-36.

%\bibitem[Re3]{Re3} \bysame, {\em Cuntz-like algebras}, Operator
%  theoretical methods, 371-386, Theta Found., Bucharest, 2000.



%\bibitem[Re3]{Re3} \bysame, {\em Cartan Subalgebras in
%    $C^*$-algebras}, Irish Math. Soc. Bulletin {\bf 61} (2008), 29-63.

%\bibitem[Re2]{Re2} \bysame, {\em  Two applications of the dual groupoid of a $C^*$-algebra}, LNM 1132, Berlin, Springer Verlag, 1985.

%\bibitem[Re3]{Re3} \bysame, {\em The Radon-Nikodym problem for approximately proper equivalence relations},  Ergodic Theory Dynam. Systems {\bf 25}(2005), 1643--1672. 

 
%\bibitem[Ri0]{Ri0} M. A. Rieffel, {\em Induced representations of $C^*$-algebras}, Adv. in Math. {\bf 13} (1974), 176-257.


%\bibitem[Ri1]{Ri1} M. A. Rieffel, {\em Strong Morita equivalence of certain transformation group $C\sp*$-algebras}, Math. Ann. {\bf 222} (1976), 7--22.


%\bibitem[Ri2]{Ri2} \bysame, {\em $C^*$-algebras associated with irrational rotations}, Pac. J. Math. {\bf 93} (1981), 415--429.

%\bibitem[Ri3]{Ri3} \bysame, {\em Non-commutative tori - a case study of non-commutative differentiable manifolds}, Contemporary Mathematics {\bf 105} (1990), 191--211.



%\bibitem[R]{R} V. A. Rohlin, {\em On the fundamental ideas of measure theories}, Translations of the American Mathematical Society (Series 1) {\bf 10} (1962), 1--54.

%\bibitem[RS]{RS} J. Rosenberg and C. Schochet, {\em The K\"unneth theorem and the universal coefficient theorem for Kasparov's generalized K-functor}, Duke J. Math. {\bf 55} (1987), 337--347.

%\bibitem[Ro]{ro} I. Rosenholtz, {\em Local expansions, derivatives, and fixed points},  Fund. Math. {\bf 91}(1976), 1--4.

%\bibitem[Ru]{Ru} S. Ruette, {\em On the Vere-Jones classifiaction and existence of maximal measures for countable topological Markov chains}, Pac. J. Math. {\bf 209} (2003), 365--380.

%\bibitem[Ru1]{Ru1} D. Ruelle, {\em Thermodynamic Formalism}, Encyclopedia of Mathematics and its Applications 5, Addison-Wesley, Reading, Ma. 1978.



%\bibitem[Ru2]{Ru2} \bysame, {\em Non-commutative algebras for hyperbolic diffeomorphisms}, Invent. Math. {\bf 93} (1988), 1--13.

\bibitem[R1]{R1} M. R\o rdam, {\em Classification of Certain Infinite
    Simple $C^*$-algebras}, J. Func. Anal. {\bf 131} (1995), 415-458. 

\bibitem[R2]{R2} \bysame, {\em Classification of Certain Infinite
    Simple $C^*$-algebras, III}, in 'Operator Algebras and Their
  Applications', Fields Institute Communications, 1997.

%\bibitem[Sa]{Sa} K.Sakai, {\em Periodic points of positively expansive maps},  Proc. Amer. Math. Soc. {\bf 94}(1985), 531--534. 

%\bibitem[S]{S} K. Schmidt, {\em Dynamical Systems of Algebraic Origin}, Birkh\"auser, Basel, Boston, Berlin, 1995.

%\bibitem[Sch]{Sch} M. Schraudner, {\em On the algebraic properties of the automorphism groups of countable-state Markov shifts}, Ergodic Theory Dynam. Systems {\bf 26} (2006), 551--583.

%\bibitem[S]{S} M. Shub, {\em Global Stability of Dynamical Systems}, Springer Verlag, Berlin, 1987.  

%\bibitem[S]{S} M. Shub, {\em Endomorphisms of compact differentiable manifolds}, Amer. J. Math. {\bf 91} (1969), 175--199.


%\bibitem[Sc]{Sc} J. Schmeling, {\em Symbolic dynamics for $\beta$-shifts and self-normal numbers}, Ergod. Th. \& Dynam. Sys. {\bf 17} (1997), 675--694.
 
%\bibitem[Sc]{Sc} M. Schraudner, {\em On the algebraic properties of the automorphism groups of countable state Markov shifts}, Preprint, Heidelberg, 2005.


\bibitem[St]{St} P.J. Stacey, {\em Crossed products of $C^*$-algebras
    by endomorphisms}, J. Austral. Math. Soc. {\bf 54} (1993), 204-212.

\bibitem[Th1]{Th1} K. Thomsen, {\em Semi-\'etale groupoids and
  applications}, Annales de l'Institute Fourier {\bf 60} (2010), 759-800.

%\bibitem[Th2]{Th2} \bysame, {\em On the $C^*$-algebra of a locally
%    injective surjection and its KMS states}, Comm. Math. Phys., to appear.

\bibitem[Th2]{Th2} \bysame,  {\em Inductive limits of interval algebras : The simple case}, in H. Araki et al (eds.), Quantum and Non-commutative Analysis (1993), 399-404.


%\bibitem[Th1]{Th1} K. Thomsen, {\em $C^*$-algebras of homoclinic and
%    heteroclinic structure in expansive dynamics}, Preprint, IMF
%  Aarhus University, 2007.

%\bibitem[Th2]{Th2} \bysame, {\em Traces, Unitary Characters and Crossed Products by $\mathbb Z$}, Publ. RIMS {\bf 31}(1995), 1011--1029.

%\bibitem[Th3]{Th3} \bysame, {\em On the KK-theory and the E-theory of
%    amalgamated free products of $C^*$-algebras}, J. Func. Anal. {\bf
%    201} (2003), 30-56.



%\bibitem[Th2]{Th2} \bysame, {\em Homomorphisms between finite direct sums of circle algebras}, Linear and Multilinear Algebra {\bf 32}(1992), 33--50.

%\bibitem[Th1]{Th1} K. Thomsen, {\em On free transformation groups and $C^*$-algebras}, Proc. Royal Soc. Edinburgh {\bf 107A} (1987), 339--347.

%\bibitem[Th2]{Th2} \bysame, {\em Inductive limits of interval algebras: the simple case},  Quantum and non-commutative analysis (Kyoto, 1992), 399--404, Math. Phys. Stud., 16, Kluwer Acad. Publ., Dordrecht, 1993.



%\bibitem[Th4]{Th4} \bysame, {\em From traces to states of the $K_0$ group of a simple $C^*$-algebras}, Bull. London Math. Soc. {\bf 28} (1996), 66--72.

%\bibitem[Th5]{Th5} \bysame, {\em On the structure of a sofic shift space}, Trans. Amer. Math. Soc. {\bf 356} (2004), 3557--3619.


%\bibitem[Th1]{Th1} K. Thomsen, {\em On the ergodic theory of synchronized systems}, Preprint, Aarhus, 2004.


%\bibitem[Th]{Th1} K. Thomsen, {\em On the structure of a sofic shift space}, Trans. Amer. Math. Soc. {\bf 356} (2004), 3557--3619.




%\bibitem[Th2]{Th2} \bysame, {\em Asymptotic Homomorphisms and Equivariant KK-theory}, J. Func. Anal. {\bf 163} (1999), 324--343.

%\bibitem[Th1]{Th1} K. Thomsen, {\em Semi \'etale groupoids and
 %   applications}, Preprint, January, 2009.

%\bibitem[Th2]{Th2} \bysame, {\em The variational principle for the defect of factor maps}, Israel J. Math. {\bf 110} (1999), 359--369.


%\bibitem[To]{To} J. Tomiyama, {\em Invitation to $C^*$-algebras and topological dynamics}, World Scientific Advanced Series in Dynamical Systems {\bf 3}, World Scientific, Singapore/New Jersey/Hong Kong, 1987.


\bibitem[T]{T} Toms, A. S., {\em Comparison theory and smooth minimal
    $C^*$-dynamics}, Comm. Math. Phys. {\bf 289} (2009), 401-433.


%\bibitem[T]{T} P. Trow, {\em Determining presentations of sofic shifts}, Theor. Comp. Science {\bf 259} (2001), 199--216.

%\bibitem[Y1]{Y1} I. Yi, {\em Canonical symbolic dynamics for
%    one-dimensional generalized solenoids},  Trans. Amer. Math. Soc.
%  {\bf 353} (2001), 3741--3767.

%\bibitem[Y2]{Y2} \bysame , {\em Ordered group invariants for
%    one-dimensional spaces},  Fund. Math. {\bf 170} (2001), 267--286. 

%\bibitem[Y3]{Y3} \bysame ,{\em Ordered group invariants for
%    nonorientable one-dimensional generalized solenoids},
%  Proc. Amer. Math. Soc. {\bf 131} (2003), 1273--1282. 

%\bibitem[Y4]{Y4} \bysame , {\em $K$-theory of $C\sp *$-algebras from
%    one-dimensional generalized solenoids},  J. Operator Theory {\bf
%    50} (2003), 283--295. 

%\bibitem[Y5]{Y5} \bysame , {\em Bratteli-Vershik systems for
%    one-dimensional generalized solenoids}, Houston J. Math. {\bf 30} (2004), 691--704.

 
%\bibitem{W1}[W1] P. Walters, {\em Relative pressure, relative equilibrium states, compensation functions and many-to-one codes between subshifts}, Trans. Amer. Math. Soc. {\bf 296} (1986), 1--31.


%\bibitem[W]{W}  P. Walters, {\em An Introduction to Ergodic Theory}, Springer Verlag, New York, Heidelberg, Berlin, 1982.

%\bibitem[W]{W} P. Walters, {\em Topological conjugacy of affine transformations of compact abelian groups},  Trans. Amer. Math. Soc. {\bf 140} (1969), 95--107.


%\bibitem[Wa]{Wa} J.B. Wagoner, {\em Topological Markov Chains, $C^*$-algebras, and $K_2$}, Advances in Math. {\bf 71} (1988), 133--185.

%\bibitem[Wi1]{Wi1} R.F. Williams, {\em One-dimensional non-wandering
%    sets}, Topology {\bf 6} (1967), 473-487.

%\bibitem[Wi2]{Wi2}  \bysame, {\em Classification of one dimensional attractors}, Proceedings of Symposia in Pure Mathematics {\bf 14} (1970), 341--361.

%\bibitem[Wi2]{Wi2} \bysame, {\em Expanding attractors}, IHES Publ. Math. {\bf 43}(1974), 169--203.

\bibitem[W]{W} W. Winter, {\em Nuclear dimension and $\mathcal
    Z$-stability of perfect $C^*$-algebras}, arXiv:1006.2731v1.


%\bibitem[Z]{Z} G. Zeller-Meyer, {\em Produits crois\'es d'une $C^*$-alg\`ebre par une group d'automorphisms}, J. Math. Pures Appl. {\bf 47}(1968), 101--239.


\end{thebibliography}
%\end{document}

\end{document}